\documentclass[english]{article}
\usepackage[T1]{fontenc}
\usepackage[latin9]{inputenc}
\usepackage{color}
\usepackage{array}
\usepackage{amsmath}
\usepackage{amssymb}
\usepackage{graphicx}

\makeatletter

\providecommand{\tabularnewline}{\\}

\usepackage{etoolbox}

\newtoggle{epstablesarenicer}
\toggletrue{epstablesarenicer}
\usepackage{colortbl}

\@ifundefined{textcolor}
 {\usepackage{color}}{}
\definecolor{gray025}{gray}{0.25}
\definecolor{gray033}{gray}{0.33}
\definecolor{gray05}{gray}{0.5}
\definecolor{gray075}{gray}{0.75}
\definecolor{gray09}{gray}{0.9}
\definecolor{gray097}{gray}{0.97}

\makeatother

\usepackage{babel}
\begin{document}

\title{Counterexamples to Cantorian Set Theory}

\author{E. Coiras}

\date{24 April 2014}
\maketitle
\begin{abstract}
This paper provides some counterexamples to Cantor's contributions
to the foundations of Set Theory. The paper starts with a very basic
counterexample to Cantor's Diagonal Method (DM) applied to binary
fractional numbers that forces it to yield one of the numbers in the
target list. To study if this specific case is just an anomaly or
rather a deeper source for concern, and given that for the DM to work
the list of numbers have to be written down, the set of numbers that
can be represented using positional fractional notation, $\mathbb{W}$,
is properly characterized. It is then shown that $\mathbb{W}$ is
not isomorphic to $\mathbb{R}$, the set of all real numbers. This
fact means that results obtained from the application of the DM to
$\mathbb{W}$ in order to derive properties of $\mathbb{R}$ are not
valid. After showing that $\mathbb{W}$ can be trivially well-ordered,
application of the DM to some shuffles of the ordered list of the
set of binary numbers $\mathbb{W_{\textrm{2}}}$ is proven to always
converge to a number in the list. This result is then used to generate
a counterexample to Cantor's DM for a generic list of reals that forces
it to yield one of the numbers of the list, thus invalidating Cantor's
result that infers the non-denumerability of $\mathbb{R}$ from the
application of the DM to $\mathbb{W}$. After this apparently anomalous
result we are forced to question Cantor's Theorem about the different
cardinalities of a set and its power set, and by means of another
counterexample we show that Cantor's Theorem does not actually hold
for infinite sets. After analyzing all these counterexamples, it is
shown that the current notion of cardinality for infinite sets does
not depend on the \textquotedbl{}size\textquotedbl{} of the sets,
but rather on the representation chosen for them. Following this line
of thought, the concept of model as a framework for the construction
of the representation of a set is introduced, and a theorem is proven
showing that an infinite set can be well-ordered if there is a proper
model for it. To reiterate that the cardinality of a set does not
determine whether the set can be well-ordered, a set of cardinality
{\normalsize{$\aleph_{0}^{\aleph_{0}}=2^{\aleph_{0}}=\mathfrak{c}$}}
is proven to be equipollent to the set of natural numbers $\mathbb{N}$.
The paper concludes with an analysis of the cardinality of the ordinal
numbers, for which a representation of cardinality {\normalsize{$\aleph_{0}^{\aleph_{0}}$}}
is proposed.
\end{abstract}

\section{Introduction}

The contributions of Georg Cantor to Set Theory \cite{Cantor1915,Zermelo1932}
have been controversial since their very inception. Perhaps his most
shocking and counterintuitive result is the apparent existence of
different types of infinity. His Diagonal Argument to show the different
cardinalities of the natural and real numbers is particularly noteworthy,
and has pervaded other areas of mathematics such as Logic and Computability
theory.

In this paper we will see that the cardinality of infinite sets seems
to be determined by the representation we choose for them, rather
than by their \textquotedbl{}size\textquotedbl{}. We will also demonstrate
that Cantor's Diagonal Method (DM) can be forced to produce one of
the listed numbers under some conditions.

Since attacking such a controversial yet widely accepted theory will
surely face some strong opposition, the paper is mostly structured
around counterexamples that are easy to explain and verify. The paper
follows the research path the author, which should make it more approachable
to general readers. Those familiar with Cantorian Set Theory may want
to skip directly to Section \ref{sec:Apply-Ldi-to-R}.

In this paper we will start with a minimalistic counterexample to
the DM and analyze it to see where it leads. We will find that the
DM works not on the set of reals $\mathbb{R}$, but on the set of
writable numbers $\mathbb{W}$. The fact that these sets are not isomorphic
is the source of the DM's apparently paradoxical consequences. In
particular it will be proved that for some orderings of the set $\mathbb{W_{\textrm{2}}}$
of numbers writable in base 2 the DM always yields a number on the
list. Using this result and the fact that the DM needs to enter $\mathbb{W}_{b}$
for some base \emph{b} will permit us to show in section \ref{sec:Apply-Ldi-to-R}
that there are shufflings of the list of all numbers in $\mathbb{R}$
that force the DM to yield one of the numbers of the list. This result
effectively invalidates the most widely accepted proof of the non-denumerability
of the real numbers.

This erratic behavior of Cantor's DM leads us to questioning the validity
of Cantor's Theorem, which states that a set and its power-set always
have different cardinality. While this is obviously true for finite
sets, Cantor relied on the DM to demonstrate its validity for infinite
sets; and after seeing the counterexamples for $\mathbb{W_{\textrm{2}}}$
and $\mathbb{W_{\textrm{10}}}$ for the special case of the non-denumerability
of the reals we are forced to examine his general result in more detail.
Ultimately we are able to provide a specific counterexample to Cantor's
Theorem by producing a bijection between the set of non-negative integers
and its power-set in section \ref{sec:Counter-Cantor-Theorem}.

In section \ref{sub:The-Applicative-Numbering} we introduce the Applicative
Numbering model to sort a set of size $\aleph_{0}^{\aleph_{0}}=2^{\aleph_{0}}=\mathfrak{c}$.
This proves that cardinality by itself is not a limiting factor for
whether a set can be well-ordered or not. In fact, we see in section
\ref{sec:Numbering-models} that there are very specific conditions
that a representation of a set must verify for it to be well-orderable.
The most important conclusion of the paper is that the standard positional
numbering model (section \ref{sub:Standard-numbering-model}) does
not meet one of the conditions, which is the main reason why the DM
has been incorrectly taken to imply that $\mathbb{R}$ is non-denumerable.

To conclude, the paper presents some tentative models for the real
and the ordinal numbers, and implications to the Continuum Hypothesis
are discussed.

\section{A minimalistic counterexample to Cantor's Diagonal Method}

Cantor claims that for every infinite list of reals in $[0,1)$, all
of them different, it is possible to find a number not contained in
the list by application of his Diagonal Method (see for instance \cite{Hodges1998}
for a basic introduction). Table \ref{tab:Table1} presents one simple
counterexample using a list $L_{1}$ of binary fractional numbers.
The third column shows the antidiagonal number $\bar{D}$ resulting
from the Diagonal Method (DM) applied to $L_{1}$ up to digit \emph{n},
where the notation $\text{�}|_{k}$ is used to denote truncation of
the string of digits up to digit \emph{k}. The fourth column shows
where in the list $L_{1}$ can the partial result $\bar{D}\left(L_{1}\right)|_{n}$
be found. 

\begin{table}[htbp]
\centering{}\iftoggle{epstablesarenicer}{\includegraphics[bb=96bp 590bp 516bp 708bp,clip,width=1\textwidth]{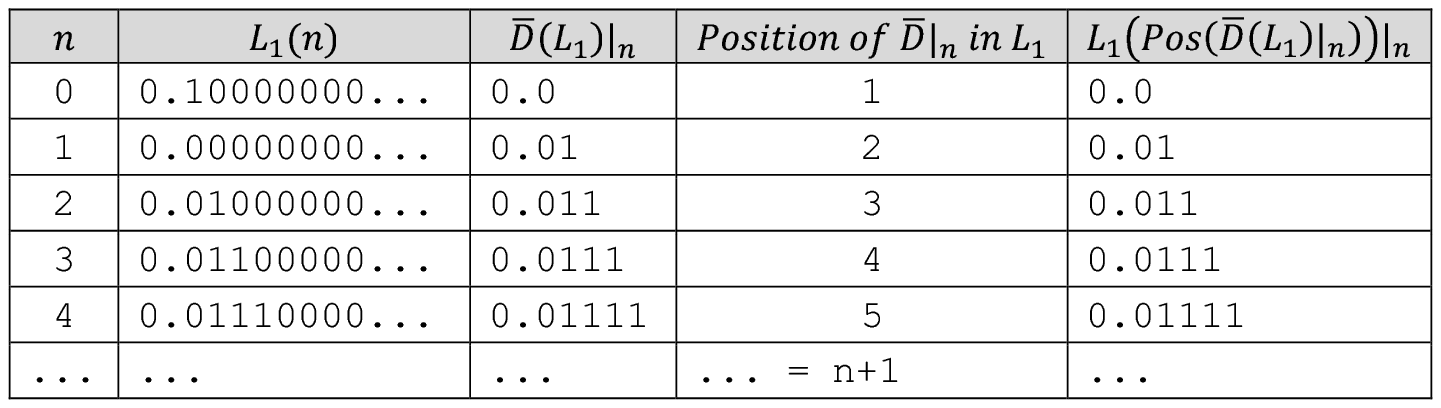}}{%
\begin{tabular}{|>{\raggedright}m{16pt}|>{\raggedright}m{66pt}|>{\raggedright}m{46pt}|>{\centering}m{58pt}|>{\raggedright}m{104pt}|}
\hline 
\rowcolor{gray09}$n$ & $L_{1}(n)$ & $\bar{D}\left(L_{1}\right)\vert_{n}$ & $Position\thinspace of$\\
$\bar{D}\vert_{n}\thinspace in\thinspace L_{1}$ & $L_{1}\left(Pos\left(\bar{D}\left(L_{1}\right)\vert_{n}\right)\right)\vert_{n}$ \tabularnewline
\hline 
0 & 0.10000000... & 0.0 & 1 & 0.0 \tabularnewline
\hline 
1 & 0.00000000... & 0.01 & 2 & 0.01 \tabularnewline
\hline 
2 & 0.01000000... & 0.011 & 3 & 0.011 \tabularnewline
\hline 
3 & 0.01100000... & 0.0111 & 4 & 0.0111 \tabularnewline
\hline 
4 & 0.01110000... & 0.01111 & 5 & 0.01111 \tabularnewline
\hline 
... & ... & ... & ... $=n+1$ & ... \tabularnewline
\hline 
\end{tabular}}\caption{A list $L_{1}(n)$ of fractional binary numbers and the partial results
of the antidiagonal number $\bar{D}$ up to digit \emph{n}. Fourth
and fifth columns show every partial antidiagonal result is already
in the list.}
\label{tab:Table1} 
\end{table}

We can see that: (1) any partial result of the DM, $\bar{D}\left(L_{1}\right)\vert_{k}$,
is equal to the next number on the list $L_{1}(k+1)$, and (2) the
limiting antidiagonal number $\bar{D}$ is equivalent to the very
first number on the list.

There is nevertheless a difference in form between the antidiagonal
number and the first number of the list: one uses 1-ending representation
and the other 0-ending. This fact was known to Cantor and Dedekind
(see \cite{Ferreiros2007}, p.191) and it's usually overlooked as
an irrelevant inconvenience to Cantor's arguments. Yet, if we are
going to be perfectly strict, we need to completely understand if
this is a cause for concern or not, and why.

For now let's just agree that the binary strings 0.10000... and 0.01111...
both represent the same real number $1/2$. After this rather trivial
observation, Table \ref{tab:Table1} seems to suggest that the DM
doesn't work as expected for at least one list $L_{1}$. And that
happens neither (1) after any finite number of steps, nor (2) after
extrapolation to the limiting result (that is, after an infinity of
steps).

Note that in Table \ref{tab:Table1} there are no dots after the partial
results, since these are finite strings of length \emph{k} for any
number of iterations \emph{k}. Note as well that these partial strings
ensure the current $\bar{D}\left(L_{1}\right)\vert_{k}$---of length
\emph{k}---is different from any of the strings found on lines 1 to
\emph{k}, that is, up to any finite number of iterations. Cantor's
argument extrapolates these partial results to the case where \emph{k}
grows to infinity, which he envisions should make sure the limiting
$\bar{D}\left(L_{1}\right)=\lim_{k\rightarrow\infty}\bar{D}\left(L_{1}\right)\vert_{k}$
is different from every string in the list, and should therefore correspond
to a completely new number not contained in the list.

However, after the admittedly simplistic counterexample shown in Table
\ref{tab:Table1}, the generalization of Cantor's extrapolation to
any list of reals doesn't seem to be as straightforward as it may
have first appeared, because, at least in principle, there is nothing
preventing every partial antidiagonal result of \emph{k} digits $\bar{D}\left(L\right)\vert_{k}$
from being found somewhere ahead of line \emph{k} in an infinite list
\emph{L}. In the example shown above this happens to occur just on
the next line, $k+1$---and for that matter on any lines after \emph{k}.
This shows the list may in fact contain a subsequence that converges
to the same limiting number $\bar{D}$. It is not clear then what
the purpose of the DM actually is: does it find one of the limits
of the list? Does the limit even belong to the same set as the numbers
being listed? These are important details that we need to study in
order to determine the soundness of Cantor's Diagonal Method.

Moreover, is the set of reals properly represented when symbolized
by a list of character strings written in positional fractional notation?
We will see in the next section that this is not the case.

\section{The set of all numbers that can be written\label{sec:The-W-set}}

The key point to understand where the problems with the DM originate
is that real numbers don't \textquotedbl{}have\textquotedbl{} digits.
For instance, $\pi$ is defined as the ratio of the length of a circle
to its diameter; there are no digits in this, in principle, purely
geometrical description. However, it is very convenient to use positional
fractional notation in base 10 to say that $\pi$ is about 3.1416,
provided we don't forget that this concatenation of digits is not
$\pi$ itself. The fact that the DM targets a list of written numbers
means that it is not targeting $\mathbb{R}$ directly, but rather
a convenient representation of it, $\mathbb{W}$. This seemingly trivial
observation has serious consequences, as we shall see. To begin with,
it means we have to be very careful when deriving properties of $\mathbb{R}$
through the analysis of the set of writable numbers $\mathbb{W}$,
especially if they happen to be completely different sets. We will
study the question within this section.

Most of what follows are well-known elementary math results. However,
I haven't been able to find them collected to specifically characterize
the set of writable numbers $\mathbb{W}$, so they are included here
for completion.

First of all, we need to define what a writable number is: a number
writable in positional fractional notation for some natural base $b>1$
is a finite string of digits smaller than \emph{b} and with a dot
separating the integer and fractional parts of the number:

\begin{equation}
\forall b>1\in\mathbb{N}:w=w_{p}\ldots w_{1}w_{0}.w_{-1}w_{-2}\ldots w_{-q}\in\mathbb{W}_{b}\Leftrightarrow p,q\in\mathbb{Z}^{+}\,\wedge\,0\le w_{k}<b\in\mathbb{Z}^{+}
\end{equation}
A $+$ or $-$ symbol may be placed at the beginning of the string
to indicate the sign of the number, but since negative numbers are
of no interest for the matter discussed in this paper, we will assume
from now on that $\mathbb{W}$ only concerns non-negative numbers.

These writable numbers are just the standard numbers we use every
day, with $b=10$ being the most common base (decimal numbers). The
numbers in base $b=2$ (binary numbers) are also important, being
the most convenient choice for computational applications.

It is easy to show that the set of all numbers that can be written
in a given base \emph{b}, $\mathbb{W}_{b}$, is not isomorphic to
the set of real numbers $\mathbb{R}$ because (1) every writable number
in $\mathbb{W}_{b}$ has an image in $\mathbb{R}$ (its value), but
(2) there are numbers in $\mathbb{R}$ that cannot be written in $\mathbb{W}_{b}$.

The real value \emph{x} of a writable number \emph{w} in $\mathbb{W}_{b}$
can be obtained by summation of the products of its digits to the
powers of \emph{b} according to the position of the digits in the
string:

\begin{equation}
\textrm{Val}:\mathbb{W}_{b}\to\mathbb{R}
\end{equation}
\begin{equation}
w=w_{p}\ldots w_{1}w_{0}.w_{-1}w_{-2}\ldots w_{-q}\in\mathbb{W}_{b}
\end{equation}
\begin{equation}
x=\textrm{Val}\left(w\right)\equiv\sum\nolimits _{k=-q}^{p}{w_{k}\cdot b^{k}}\in\mathbb{R}\label{eq:ValDef}
\end{equation}
The real number \emph{x} and its fractional positional number notation
\emph{w} are different things indeed, but this is easy to forget given
our routine use of $\mathbb{W}_{10}$ to speak about numbers in $\mathbb{R}$.
For instance, given that the string $w="123.5"$ and its corresponding
real value $x=\mathrm{Val}\left("123.5"\right)=1\cdot100+2\cdot10+3+5/10=123.5$
happen to be the same written string may make us unconsciously identify
\emph{w} with \emph{x}, and therefore $\mathbb{W}_{10}$ with $\mathbb{R}$
for any number. Yet, the number \emph{x} in $\mathbb{R}$ described
by 123.5 in decimal can also be written down as 1111011.1 in binary;
\emph{x} is neither of those strings. These seemingly trivial comments
will be very relevant when we realize that some of Cantor's derivations
suffer from an implicit association of $\mathbb{W}$ with $\mathbb{R}$.
This association breaks down when we consider numbers with an infinite
tail of digits, with the exception of the trivial terminations ...000...
and ...(\emph{b}-1) (\emph{b}-1) (\emph{b}-1)..., which do have equivalent
writable representations within $\mathbb{W}_{b}$.

It is easy to show that $\mathbb{W}_{b}$ is in fact a subset of $\mathbb{Q}$,
the set of rational numbers. The proof is a well known result:

\begin{equation}
\forall w=w_{p}\ldots w_{1}w_{0}.w_{-1}\ldots w_{-q}\in\mathbb{W}_{b}:\exists z\in\mathbb{Z},\exists n\in\mathbb{N}\vert\mathrm{Val}\left(w\right)=\frac{z}{n}\in\mathbb{Q}
\end{equation}

\[
Proof:\:\mathrm{Val}\left(w\right)=\sum\limits _{k=-q}^{p}{w_{k}\cdot b^{k}}=\frac{\sum\nolimits _{k=-q}^{p}{w_{k}\cdot b^{q+k}}}{b^{q}}=\frac{z}{n}\in\mathbb{Q}
\]

To show that $\mathbb{W}_{b}$ is a proper subset of $\mathbb{Q}$
we just need to prove that some numbers in $\mathbb{Q}$ (and therefore
in $\mathbb{R}$) cannot be represented with a finite number of digits:

\begin{equation}
\forall b>1\in\mathbb{N}\:\exists x\in\mathbb{Q}\:\vert\:\nexists w\in\mathbb{W}_{b}:\mathrm{Val}\left(w\right)=x
\end{equation}

\[
Proof:\: Take\thinspace x=\frac{1}{p}\in\mathbb{Q}\thinspace with\thinspace p\thinspace prime\thinspace and\thinspace p>b
\]
\[
Assume\thinspace\exists w\in\mathbb{W}_{b}\:\vert\: x=\mathrm{Val}\left(w\right)=\mathrm{Val}(0.w_{1}w_{2}\ldots w_{n})\thinspace with\thinspace n\in\mathbb{N}\Rightarrow
\]
\[
\Rightarrow\frac{1}{p}=\frac{w_{1}b^{n-1}+w_{2}b^{n-2}+\cdots+w_{n}}{b^{n}}\Rightarrow
\]
\[
\Rightarrow\frac{b^{n}}{p}=w_{1}b^{n-1}+w_{2}b^{n-2}+\cdots+w_{n}\in\mathbb{N}\Rightarrow p\thinspace divides\thinspace b\thinspace(contradiction)
\]
Therefore $\mathbb{W}_{b}$ is a proper subset of $\mathbb{Q}$ and
thus of $\mathbb{R}$:

\begin{equation}
\mathrm{Val}\left(\mathbb{W}_{b}\right)\subset\mathbb{Q}\subset\mathbb{R}
\end{equation}
However, the closure of $\mathbb{W}_{b}$ for any \emph{b} is the
set of all real numbers written in that base $\mathbb{R}_{b}$: 
\begin{equation}
\forall b>1\in\mathbb{N}:\mathrm{Cl}\left(\mathbb{W}_{b}\right)=\mathbb{R}_{b}
\end{equation}
And therefore the value of the closure is the set of all real numbers
$\mathbb{R}$:

\begin{equation}
\forall b>1\in\mathbb{N}:\mathrm{Val}\left(\mathrm{Cl}\left(\mathbb{W}_{b}\right)\right)=\mathrm{Val}\left(\mathbb{R}_{b}\right)=\mathbb{R}\label{eq:Cl_W_is_R}
\end{equation}
We can prove \ref{eq:Cl_W_is_R} by showing that any real number can
be approximated to any level of precision using a written positional
expansion: 
\begin{equation}
\forall x\in\mathbb{R\:}\forall\varepsilon>0\in\mathbb{R\:}\exists w\in\mathbb{W}_{b}\,|\,\left\Vert x-\mathrm{Val}(w)\right\Vert <\varepsilon\label{eq:Proof_Cl_W_is_R}
\end{equation}
\[
Proof:\thinspace\forall\varepsilon>0\:\exists k\in\mathbb{Z}^{+}\vert-k<\log_{b}\varepsilon\Rightarrow b^{-k}<\varepsilon
\]
\[
w=x\vert_{k}=x_{p}\ldots x_{1}x_{0}.x_{-1}x_{-2}\ldots x_{-k}\in\mathbb{W}_{b}\Rightarrow
\]
\[
\Rightarrow\left\Vert x-\mathrm{Val}(w)\right\Vert =0.00\ldots00x_{-k-1}x_{-k-2}\ldots<0.00\ldots01=b^{-k}<\varepsilon
\]

To complete the characterization of $\mathbb{W}_{b}$, and because
it will be of use in following sections, it is interesting to note
that the set can be expressed as the internal direct sum of the subset
of numbers with null fractional part $\mathbb{W}_{b}^{I}$ and the
subset of numbers with null integer part $\mathbb{W}_{b}^{F}$: 
\begin{equation}
\mathbb{W}_{b}=\mathbb{W}_{b}^{I}\oplus\mathbb{W}_{b}^{F}
\end{equation}
\begin{equation}
\mathbb{W}_{b}^{I}\equiv\left\{ w=w_{p}\ldots w_{1}w_{0}.w_{-1}w_{-2}\ldots w_{-q}\in\mathbb{W}_{b}\,\vert\,\mathrm{Val}\left(0.w_{-1}\ldots w_{-q}\right)=0\right\} 
\end{equation}
\begin{equation}
\mathbb{W}_{b}^{F}\equiv\left\{ w=w_{p}\ldots w_{1}w_{0}.w_{-1}w_{-2}\ldots w_{-q}\in\mathbb{W}_{b}\,\vert\,\mathrm{Val}\left(w_{p}\ldots w_{1}w_{0}.0\right)=0\right\} 
\end{equation}
Which is straightforward since every element w in $\mathbb{W}_{b}$
can be expressed as the sum of one element from $\mathbb{W}_{b}^{I}$
(its integral part) and another from $\mathbb{W}_{b}^{F}$ (its fractional
part):

\begin{equation}
w=w_{p}\ldots w_{1}w_{0}.w_{-1}w_{-2}\ldots w_{-q}=w_{p}\ldots w_{1}w_{0}+0.w_{-1}w_{-2}\ldots w_{-q}=\left\lfloor w\right\rfloor +\left\{ w\right\} 
\end{equation}

It is easy to see that the set of values of the elements in $\mathbb{W}_{b}^{I}$
is $\mathbb{Z}^{+}$, so we could write that $\mathbb{W}_{b}^{I}$
is just the representation of $\mathbb{Z}^{+}$ in base b, $\mathbb{Z}_{b}^{+}$:

\begin{equation}
\mathbb{W}_{b}^{I}=\mathrm{Base}{}_{b}\left(\mathbb{Z}^{+}\right)\equiv\mathbb{Z}_{b}^{+}
\end{equation}

\begin{equation}
\mathbb{Z}^{+}=\mathrm{Val}\left(\mathbb{W}_{b}^{I}\right)
\end{equation}
On the other hand $\mathbb{W}_{b}^{F}$ is the set of fractional real
numbers in {[}0, 1) which can be written down using finitely many
significant digits. We can denote this interval as $\left[0,1\right)_{\mathbb{W}_{b}}$
to distinguish it from the representation in base \emph{b} of the
fractional reals $\left[0,1\right)_{b}$, which contains numbers with
infinitely many significant digits:

\begin{equation}
\mathbb{W}_{b}^{F}=[0,\thinspace1)\cap\mathbb{W}_{b}\equiv[0,\thinspace1)_{\mathbb{W}_{b}}
\end{equation}
\begin{equation}
\mathrm{Val}\left([0,\thinspace1)_{\mathbb{W}_{b}}\right)\subset\mathrm{Val}\left([0,\thinspace1)_{b}\right)=[0,\thinspace1)\subset\mathbb{R}
\end{equation}
Therefore, using (11), $\mathbb{W}_{b}$ can be seen as the Cartesian
product of the set of integers in base b with the set of those fractional
numbers in {[}0, 1) that can actually be written down. This means
$\mathbb{W}_{b}$ is missing some elements present in the positional
fractional representation of $\mathbb{R}$ in base \emph{b}, $\mathbb{R}_{b}$:

\begin{equation}
\mathbb{W}_{b}=\mathbb{Z}_{b}^{+}\times[0,\thinspace1)_{\mathbb{W}_{b}}\subset\mathbb{Z}_{b}^{+}\times[0,\thinspace1)_{b}=\mathbb{R}_{b}^{+}
\end{equation}
There is nothing really new on the results shown so far in this section:
it is well known that some real numbers require a never-ending expansion
of digits when written in fractional positional notation. Some of
these are irrational or transcendental, but many others are just rationals.
Nevertheless, it was important to characterize $\mathbb{W}$ properly
in order to show more interesting results in following sections.

The differences between $\mathbb{W}$ and $\mathbb{R}$ directly mean
that it is impossible to write a list of all real numbers using positional
fractional notation, since some numbers just can't be written with
finitely many digits. Cantor could have stopped there, but instead
moved further ahead, applying the DM to a list of written numbers
without realizing he wasn't targeting $\mathbb{R}$, but a different
set, $\mathbb{W}$. He concludes then that there are ``more'' numbers
in $\mathbb{R}$ than those on any list of writable numbers (which
in a way is true, as we have seen). But from there (or, more precisely,
after he determined the cardinality of $\mathbb{R}$ to be $2^{\aleph_{0}}$
plus the application of the DM to support Cantor's Theorem) he incorrectly
infers that the reason had to be that $\mathbb{R}$ is \textquotedbl{}unlistable\textquotedbl{}
(non-denumerable). We will see in more detail why this inference is
not justified in the following sections of the paper and provide some
supporting counterexamples.

Regarding the previous paragraph, an important remark must be made:
some readers may argue that by considering tails of infinitely many
digits we can in fact target $\mathbb{R}$ instead of $\mathbb{W}$;
unfortunately this is irrelevant, because the DM only considers digits
at finite positions, so those numbers with infinite tails are completely
indistinguishable to the DM from numbers in $\mathbb{W}$. Therefore,
the only quality that could differentiate those ``unwrittable''
reals from members of $\mathbb{W}$ is not used at all by the DM.
We will later exploit this weakness of the DM to produce a counterexample
that invalidates its application to proving the non-denumerability
of $\mathbb{R}$.

We will see in the next section that $\mathbb{W}$ can actually be
well-ordered and is therefore denumerable. After that we will apply
the DM to one ordered list of $\mathbb{W}$, and find the key result
that for $\mathbb{W}_{2}$ the DM can be forced to always yield a
number explicitly written in the initial list (this is not always
true when using other bases \emph{b}>2). We will combine this result
with the weakness mentioned in the previous paragraph to construct
the equivalent counterexample for the list of all real numbers.

\section{A well-ordering for all binary writable numbers $\mathbb{W}_{2}$}

From now on I will be moving to binary numbers since they produce
more compact tables and are enough to prove that the DM cannot be
used to justify the non-denumerability of the reals. Most of the results
we'll find can be applied to other bases \emph{b}>2. Some versions
of the DM for \emph{b}=10 are studied in Section \ref{sec:Other-DM-flavours}.

To well-order $\mathbb{W}_{2}$ we first need to well order is fractional
part $\mathbb{W}_{2}^{F}=[0,\thinspace1)_{\mathbb{W}_{2}}$ by finding
a bijection that maps any non-negative integer to a number in $\mathbb{W}_{2}^{F}$
and vice versa. Fortunately this is not difficult to do: just revert
the order of the digits of every binary non-negative integer $n=x_{1}x_{2}\ldots x_{k}$
to create a fractional number between zero and one. This Digital Inversion
(DI) process operates as follows:

\begin{equation}
\mathrm{DI}\thinspace:\thinspace\mathbb{Z}^{+}\to[0,\thinspace1)_{\mathbb{W}_{2}}\subset\mathbb{W}_{2}\subset\mathbb{R}
\end{equation}
\begin{equation}
n=x_{1}x_{2}\ldots x_{k}\in\mathbb{Z}^{+},\thinspace x_{i}\in\{0,\thinspace1\}
\end{equation}
\begin{equation}
\mathrm{DI}\left(n\right)=\mathrm{DI}\left(x_{1}x_{2}\ldots x_{k}\right)\equiv0.x_{k}\ldots x_{2}x_{1}\in[0,\thinspace1)_{\mathbb{W}_{2}}
\end{equation}
An algorithm that applies the DI function to an enumeration of the
set of non-negative integers $\mathbb{Z}^{+}$ using base 2 results
in the list $L_{\mathrm{DI}}$ shown in Table \ref{tab:Table2}.

\begin{table}[htbp]
\centering{}\iftoggle{epstablesarenicer}{\includegraphics[bb=125bp 535bp 315bp 708bp,clip]{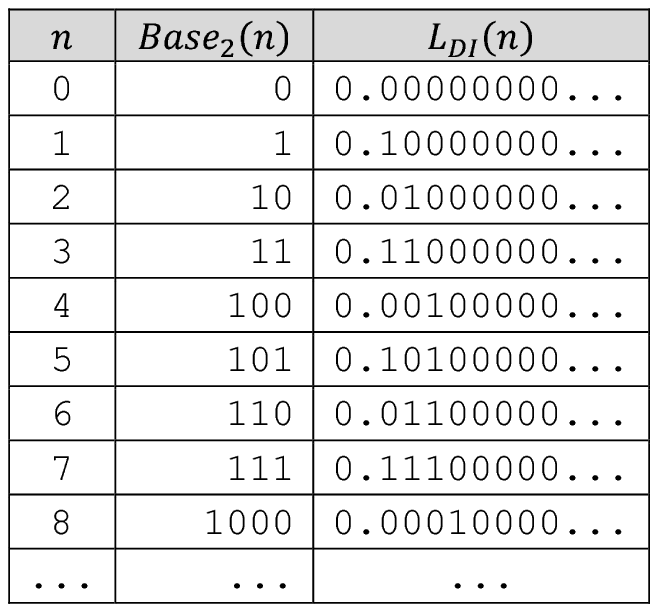}}{\begin{center}
\begin{tabular}{|>{\centering}p{30pt}|r|l|}
\hline 
\rowcolor{gray09}$n$ & $\mathrm{Base_{2}}\left(n\right)$ & $L_{\mathrm{DI}}\left(n\right)$\tabularnewline
\hline 
0 & 0 & 0.00000000... \tabularnewline
\hline 
1 & 1 & 0.10000000... \tabularnewline
\hline 
2 & 10 & 0.01000000... \tabularnewline
\hline 
3 & 11 & 0.11000000... \tabularnewline
\hline 
4 & 100 & 0.00100000... \tabularnewline
\hline 
5 & 101 & 0.10100000... \tabularnewline
\hline 
6 & 110 & 0.01100000... \tabularnewline
\hline 
7 & 111 & 0.11100000... \tabularnewline
\hline 
8 & 1000 & 0.00010000... \tabularnewline
\hline 
... & ... & ... \tabularnewline
\hline 
\end{tabular}
\par\end{center}}\caption{The list $L_{\mathrm{DI}}$ of all possible fractional numbers that
can be written in base 2.}
\label{tab:Table2} 
\end{table}

Note that for fractional numbers with finitely many significant digits
the DI function is invertible, and can therefore be used to find the
position in the list of the fractional number:

\begin{equation}
\mathrm{Pos}\left(x=0.x_{1}x_{2}\ldots x_{k}\right)=\mathrm{DI}^{-1}\left(0.x_{1}x_{2}\ldots x_{k}\right)=x_{k}\ldots x_{2}x_{1}=\thinspace n\in\mathbb{Z}^{+}
\end{equation}
The function DI therefore establishes a bijection between $\mathbb{Z}^{+}$
and $\left[0,1\right)_{\mathbb{W}_{2}}$ which suffices to produce
a well-ordering of $\left[0,1\right)_{\mathbb{W}_{2}}$.

If a fractional number has infinitely many significant digits the
DI function can still be used to find the position in $L_{\mathrm{DI}}$
of every approximation up to any digit \emph{k}:

\begin{equation}
\mathrm{DI}^{-1}\left(0.x_{1}x_{2}x_{3}\ldots\vert_{k}\right)=x_{k}\ldots x_{2}x_{1}=\thinspace n\in\mathbb{Z}^{+}
\end{equation}
Note that the way $L_{\mathrm{DI}}$ is constructed trivially guarantees
that every block of $2^{n}$ consecutive lines contains all possible
fractional numbers that can be written using any combination of n
digits. This by the way implies that any partial fractional number
$0.x_{1}x_{2}\ldots x_{k}$ can be found not only at position $n=x_{k}\ldots x_{2}x_{1}$
but also at infinitely many positions $1x_{k}\ldots x_{2}x_{1}$,
$10x_{k}\ldots x_{2}x_{1}$, $11x_{k}\ldots x_{2}x_{1}$, $100x_{k}\ldots x_{2}x_{1}$,
... ahead of line \emph{n} in $L_{\mathrm{DI}}$.

Note also that every possible real number \emph{x} in the closed interval
{[}0, 1{]} is either in $L_{\mathrm{DI}}$ or is a limit of it. The
proof is essentially the same used to show that $\mathrm{Cl}\left(\mathbb{W}_{b}\right)=\mathbb{R}$
in the previous section \ref{eq:Proof_Cl_W_is_R}.

Moving on, we just showed that $\left[0,1\right)_{\mathbb{W}_{2}}$
is denumerable, and we saw before that $\mathbb{W}_{2}$ is equivalent
to the Cartesian product of $\mathbb{Z}_{2}^{+}$ and $\left[0,1\right)_{\mathbb{W}_{2}}$.
Therefore, a list of all binary writable numbers can be composed using
Cantor's procedure for enumerating the rational numbers. In this way
we can create a 2D matrix listing on one axis the set of non-negative
integers and on the other the ordered set of all writable fractional
expansions $L_{\mathrm{DI}}$as shown in Table \ref{tab:Table3}.

\begin{table}[htbp]
\centering{}\iftoggle{epstablesarenicer}{\includegraphics[bb=78bp 567bp 495bp 708bp,clip,width=1\textwidth]{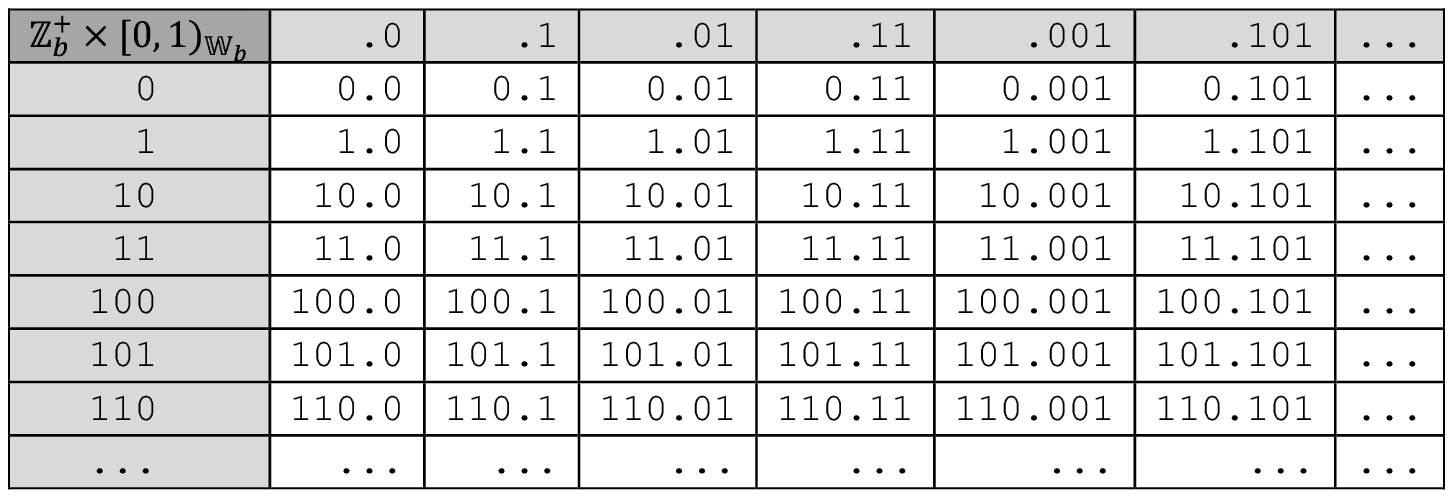}}{%
\begin{tabular}{|>{\raggedleft}p{23pt}|r|r|r|r|r|r|r|}
\cline{2-8} 
\multicolumn{1}{>{\raggedleft}p{23pt}|}{} & \cellcolor{gray09}.0 & \cellcolor{gray09}.1 & \cellcolor{gray09}.01 & \cellcolor{gray09}.11 & \cellcolor{gray09}.001 & \cellcolor{gray09}.101 & \cellcolor{gray09}... \tabularnewline
\hline 
\cellcolor{gray09}0 & 0.0 & 0.1 & 0.01 & 0.11 & 0.001 & 0.101 & ... \tabularnewline
\hline 
\cellcolor{gray09}1 & 1.0 & 1.1 & 1.01 & 1.11 & 1.001 & 1.101 & ... \tabularnewline
\hline 
\cellcolor{gray09}10 & 10.0 & 10.1 & 10.01 & 10.11 & 10.001 & 10.101 & ... \tabularnewline
\hline 
\cellcolor{gray09}11 & 11.0 & 11.1 & 11.01 & 11.11 & 11.001 & 11.101 & ... \tabularnewline
\hline 
\cellcolor{gray09}100 & 100.0 & 100.1 & 100.01 & 100.11 & 100.001 & 100.101 & ... \tabularnewline
\hline 
\cellcolor{gray09}101 & 101.0 & 101.1 & 101.01 & 101.11 & 101.001 & 101.101 & ... \tabularnewline
\hline 
\cellcolor{gray09}110 & 110.0 & 110.1 & 110.01 & 110.11 & 110.001 & 110.101 & ... \tabularnewline
\hline 
\cellcolor{gray09}... & ... & ... & ... & ... & ... & ... & ... \tabularnewline
\hline 
\end{tabular}}\caption{The Cartesian product of the integer and fractional parts of $\mathbb{W}_{2}$.}
\label{tab:Table3} 
\end{table}

We then traverse it diagonally to obtain a sorted list of all binary
writable numbers in $\mathbb{W}_{2}$: 
\begin{equation}
L\left(\mathbb{W}_{2}\right)=\left\{ 0.0,\thinspace0.1,\thinspace1.0,\thinspace0.01,\thinspace1.1,\thinspace10.0,\thinspace0.11,\thinspace1.01,\thinspace10.1,\thinspace11.0,\thinspace0.001,\thinspace1.11,\thinspace10.01,\thinspace11.1,\thinspace100.0,\thinspace\ldots\right\} 
\end{equation}
Since $L\left(\mathbb{W}_{2}\right)$ is just an ordering of the elements
of $\mathbb{W}_{2}$, the closure of $L\left(\mathbb{W}_{2}\right)$
is the set of all real binary numbers $\mathbb{R}_{2}$, as we saw
in section \ref{sec:The-W-set}: 
\begin{equation}
\mathrm{Cl}\left(L\left(\mathbb{W}_{2}\right)\right)=\mathrm{Cl}\left(\mathbb{W}_{2}\right)=\mathbb{R}_{2}
\end{equation}
These results can be easily generalized to any other base different
than 2, meaning that all $\mathbb{W}_{b}$ can be well-ordered and
are therefore denumerable.

Note however that the closure of $\mathbb{W}_{b}$ cannot be well-ordered
in the same lexicographic fashion as $\mathbb{W}_{b}$ is, since infinitely
many elements of the closure require infinitely many significant digits
and therefore correspond to the same limiting \textquotedbl{}position\textquotedbl{}
of infinity (this is what Cantor actually found, but that the reals
are not denumerable does not follow from this fact, as we shall see
in section \ref{sub:Standard-numbering-model}). Nevertheless, some
numbers in the closure's boundary are equivalent in $\mathbb{R}$
(that is, have the same value) to members of $\mathbb{W}_{b}$, namely
all those that end in an infinite tail of repeated 0 digits or of
repeated (\emph{b}-1) digits. This will be essential to demonstrate
that the DM fails in general to provide a non-listed number in base
2, as presented in the next section.

I'll be using the notation $L\left(S\right)$ throughout the paper
to denote a listing of a given set \emph{S}, as done above without
introduction. This will imply that \emph{S} admits a well-ordering,
that is, can be put into one-to-one correspondence with the set of
non-negative integers. We may then write $L_{\mathrm{DI}}=L\left([0,\thinspace1)_{\mathbb{W}_{2}};\,\mathrm{DI}\right)$
where the element after the semicolon indicates the mapping function.
We can make the well-ordering property for a set S explicit by writing
$\exists L\left(S\right)$.

\section{Application of the DM to $L_{\mathrm{DI}}$}

Direct application of the DM to $L_{\mathrm{DI}}$ results in the
limiting number $\overline{D}=0.11111\ldots$ which doesn't belong
to {[}0, 1) nor is explicitly written on the list, even if every partial
expansion $0.11111\ldots|{}_{\mathrm{k}}$ is. A simple way to make
it converge to a number in a finite position of the list is to apply
a shuffle $S$ to the rows of $L_{\mathrm{DI}}$, generating a new
list $L_{\mathrm{DI}}^{'}$. Consider this very simple shuffle $S_{0}$:
\begin{equation}
S_{0}\thinspace:\thinspace\mathbb{Z}^{+}\to\mathbb{Z}^{+}
\end{equation}
\begin{equation}
S_{0}\left(n\right)=\left\{ \begin{array}{cc}
1-n & if\thinspace n\le1\\
n & if\thinspace n>1
\end{array}\right.
\end{equation}
When applied to $L_{\mathrm{DI}}$, $S_{0}$ swaps the first two rows.
Let's define $L_{\mathrm{DI}}^{'}\left(n\right)=\left(S_{0}\ast L_{\mathrm{DI}}\right)\left(n\right)=L_{\mathrm{DI}}\left(S_{0}\left(n\right)\right)$
and see what the application of the DM yields for this case; the result
is shown in Table \ref{tab:Table4}.

\begin{table}[htbp]
\centering{}\iftoggle{epstablesarenicer}{\includegraphics[bb=81bp 424bp 506bp 708bp,clip,width=1\textwidth]{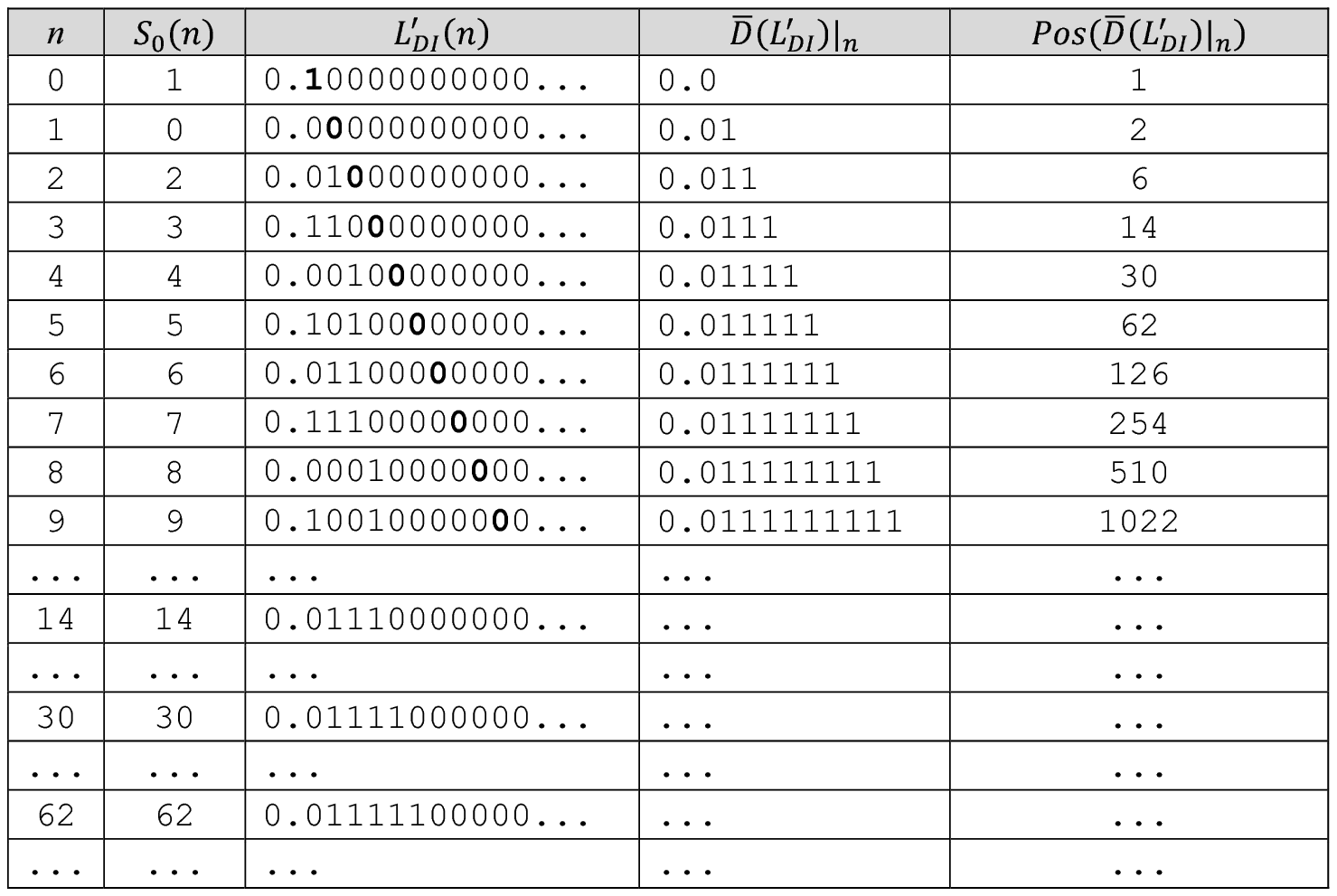}}{%
\begin{tabular}{|>{\centering}p{20pt}|c|l|>{\raggedright}p{66pt}|>{\centering}p{80pt}|}
\hline 
\rowcolor{gray09}$n$ & $S_{0}\left(n\right)$ & $L_{\mathrm{DI}}^{'}\left(n\right)$ & $\bar{D}\left(L_{\mathrm{DI}}^{'}\right)\vert_{n}$ & $\mathrm{Pos}\left(\bar{D}\left(L_{\mathrm{DI}}^{'}\right)\vert_{n}\right)$ \tabularnewline
\hline 
0 & 1 & 0.\textbf{1}0000000000... & 0.0 & 1 \tabularnewline
\hline 
1 & 0 & 0.0\textbf{0}000000000... & 0.01 & 2 \tabularnewline
\hline 
2 & 2 & 0.01\textbf{0}00000000... & 0.011 & 6 \tabularnewline
\hline 
3 & 3 & 0.110\textbf{0}0000000... & 0.0111 & 14 \tabularnewline
\hline 
4 & 4 & 0.0010\textbf{0}000000... & 0.01111 & 30 \tabularnewline
\hline 
5 & 5 & 0.10100\textbf{0}00000... & 0.011111 & 62 \tabularnewline
\hline 
6 & 6 & 0.011000\textbf{0}0000... & 0.0111111 & 126 \tabularnewline
\hline 
7 & 7 & 0.1110000\textbf{0}000... & 0.01111111 & 254 \tabularnewline
\hline 
8 & 8 & 0.00010000\textbf{0}00... & 0.011111111 & 510 \tabularnewline
\hline 
9 & 9 & 0.100100000\textbf{0}0... & 0.0111111111 & 1022 \tabularnewline
\hline 
... & ... & ... & ... & ... \tabularnewline
\hline 
14 & 14 & 0.01110000000... & ... & ... \tabularnewline
\hline 
... & ... & ... & ... & ... \tabularnewline
\hline 
30 & 30 & 0.01111000000... & ... & ... \tabularnewline
\hline 
... & ... & ... & ... & ... \tabularnewline
\hline 
62 & 62 & 0.01111100000... & ... & ... \tabularnewline
\hline 
... & ... & ... & ... & ... \tabularnewline
\hline 
\end{tabular}}\caption{A shuffling of the list of all writable numbers in {[}0, 1) and the
result of the application of the DM.}
\label{tab:Table4} 
\end{table}

And we see that the value of the limiting antidiagonal number $\bar{D}$
actually matches that of the first number on the list; also, any partial
result of the DM can be found somewhere ahead on the list. Here, again,
we need to be especially careful in keeping in mind that the strings
in $\mathbb{W}_{2}$ and the numbers they represent in $\mathbb{R}$
are different objects. So, it doesn't matter that $L_{\mathrm{DI}}^{'}\left(0\right)=0.10000\ldots$
and $\bar{D}\left(L_{\mathrm{DI}}^{'}\right)=0.01111\ldots$ are different
binary strings, because their values in $\mathbb{R}$ are the same:
$\mathrm{Val}\left(L_{\mathrm{DI}}^{'}\left(0\right)\right)=\mathrm{Val}\left(\bar{D}\left(L_{\mathrm{DI}}^{'}\right)\right)=1/2\in\mathbb{R}$

\begin{figure}
\begin{centering}
\includegraphics[bb=65bp 517bp 550bp 700bp,clip,scale=0.55]{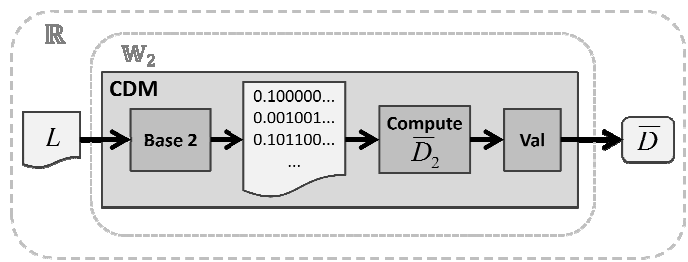}
\par\end{centering}

\caption{A block diagram representation of Cantor's Diagonal Method (CDM) using
base 2.}
\label{fig:Figure1}
\end{figure}

To make this distinction clear we can think of the DM as a black box
that takes a list of reals $L$ and outputs a real number $\bar{D}$
(Figure \ref{fig:Figure1}). Inside the box we first have to convert
the list of reals to a list of strings in $\mathbb{W}_{b}$ using
a base \emph{b} positional fractional notation; from those strings
we obtain the diagonal number $D_{b}$, and from it the antidiagonal
number in base \emph{b}, $\overline{D}_{b}$. Application of the Val
function to $\overline{D}_{b}$ yields the final output $\bar{D}\in\mathbb{R}$.

So, to summarize what we saw in this section, we started with an infinite
list $L_{\mathrm{DI}}^{'}$ of reals, all of them different, put them
through a black box called DM and obtained an output $\bar{D}$ whose
value actually matches that of one of the reals in $L_{\mathrm{DI}}^{'}$.
This seems to indicate that Cantor's Diagonal Argument can't possibly
work in general, since we have found at least one specific case where
it fails to produce a number not on the input list. Furthermore, $L_{\mathrm{DI}}^{'}$
contains every fractional number that can be written, so it is not
at all clear how the situation could be more or less favorable for
the DM if we consider as well the missing unwritable reals in {[}0,
1).

It may be that this counterexample is just an anomaly, and that the
DM actually works well for the list of all reals $L\left(\mathbb{R}\right)$.
But we will see in section \ref{sec:Apply-Ldi-to-R} that this is
not the case: the DM also fails for the list of all reals in {[}0,
1). Before that, however, we need to address a minor detail on the
implementation of the DM black box.

\section{Addressing other variants of Cantor's Diagonal Method\label{sec:Other-DM-flavours}}

Some readers may point out that the lists provided so far are \textquotedbl{}wrong\textquotedbl{},
because we are supposed to make sure the numbers on the list use the
1-ending representation of the real numbers that have the infinite
termination ...\emph{xyz}01111... instead of the equivalent 0-ending
representation ...\emph{xyz}10000... This requirement appears in many
texts discussing the DM (e.g. \cite{Huntington2003,Hawking2006})
and seems to have originated in 1877 when Dedekind pointed it out
to Cantor as a potential problem with Cantor's proof of the equipollence
of $\mathbb{R}$ and $\mathbb{R}^{n}$ (see \cite{Ferreiros2007},
p.191).

If we are to use the 1-ending criterion we may have to discard 0 and
target (0, 1) or (0, 1{]} instead of {[}0, 1), because 0 doesn't have
an equivalent 1-ending representation. In any case, applying any of
these extra conditions doesn't make any significant difference to
the performance of the DM. To prove it let's create a version of the
list without 0.0 in it and using the 1-ending version of the numbers
to obtain a new list $L_{\mathrm{DI}}^{''}$ meeting these more stringent
requirements.

The result is presented in Table \ref{tab:Table5} where it can be
seen that the DM also yields a listed number when using the 1-ending
representation for the input list. Note that in this case we didn't
even need to add a shuffle function to make it converge to a number
explicitly on the list. Here I must again remind the reader that the
strings in $L_{\mathrm{DI}}^{''}$ and $\bar{D}$ are in positional
notation and we still need to apply the Val function (\ref{eq:ValDef})
to obtain the actual numbers they represent in $\mathbb{R}$. In this
way we observe $\mathrm{Val}\left(L_{\mathrm{DI}}^{''}\left(3\right)\right)=\mathrm{Val}\left(\bar{D}\left(L_{\mathrm{DI}}^{''}\right)\right)=3/4\in\mathbb{R}$.

\begin{table}[htbp]
\centering{}\iftoggle{epstablesarenicer}{\includegraphics[bb=98bp 537bp 429bp 708bp,clip]{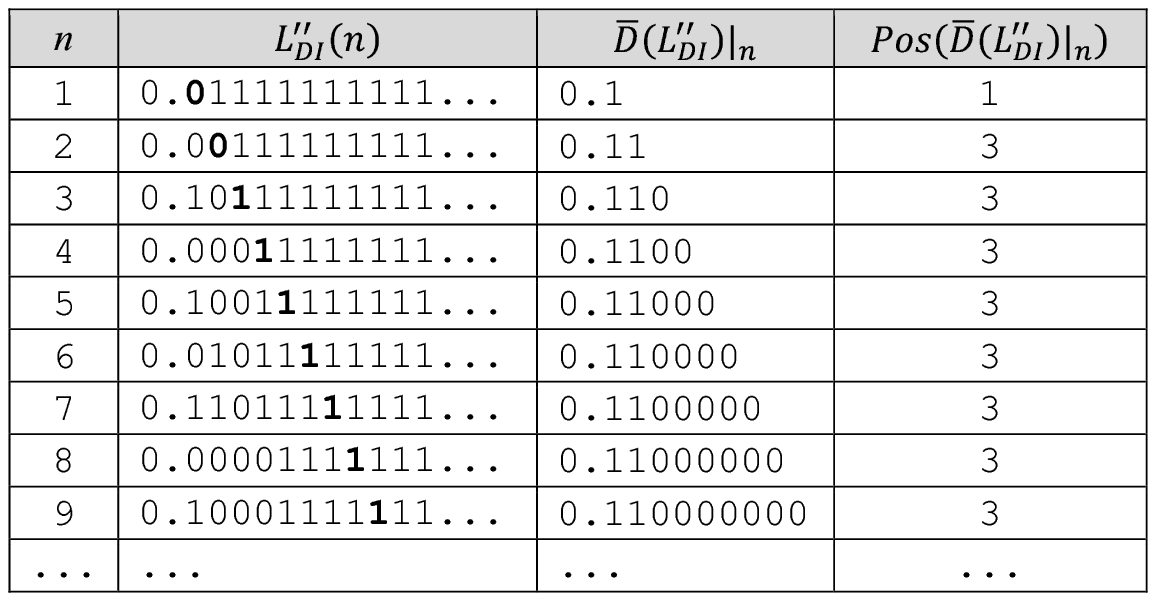}}{%
\begin{tabular}{|>{\centering}p{25pt}|>{\raggedright}p{86pt}|>{\raggedright}p{66pt}|>{\centering}p{82pt}|}
\hline 
\rowcolor{gray09}$n$ & $L_{\mathrm{DI}}^{''}\left(n\right)$ & $\bar{D}\left(L_{\mathrm{DI}}^{''}\right)\vert_{n}$ & $\mathrm{Pos}\left(\bar{D}\left(L_{\mathrm{DI}}^{''}\right)\vert_{n}\right)$ \tabularnewline
\hline 
1 & 0.\textbf{0}1111111111... & 0.1 & 1 \tabularnewline
\hline 
2 & 0.0\textbf{0}111111111... & 0.11 & 3 \tabularnewline
\hline 
3 & 0.10\textbf{1}11111111... & 0.110 & 3 \tabularnewline
\hline 
4 & 0.000\textbf{1}1111111... & 0.1100 & 3 \tabularnewline
\hline 
5 & 0.1001\textbf{1}111111... & 0.11000 & 3 \tabularnewline
\hline 
6 & 0.01011\textbf{1}11111... & 0.110000 & 3 \tabularnewline
\hline 
7 & 0.110111\textbf{1}1111... & 0.1100000 & 3 \tabularnewline
\hline 
8 & 0.0000111\textbf{1}111... & 0.11000000 & 3 \tabularnewline
\hline 
9 & 0.10001111\textbf{1}11... & 0.110000000 & 3 \tabularnewline
\hline 
... & ... & ... & ... \tabularnewline
\hline 
\end{tabular}}\caption{A binary representation of the list of all writable numbers in {[}0,
1) using 1-ending termination and the result of the application of
Cantor's Diagonal Method.}
\label{tab:Table5} 
\end{table}

Other flavors of the DM concern the number base used. In base 10,
which is the variant most commonly used in popular math books, you
can prevent the antidiagonal number from converging to an explicit
number on the list by choosing the replacement digits properly. Table
\ref{tab:Table6} summarizes the results for some implementations
of the DM in base 10 found in math vulgarization books when applied
to $L_{\mathrm{DI}}$. The most commonly found version is the one
that adds 1 to every digit modulo 10, and can be found in Stephen
Hawking's \textquotedbl{}God Created the Integers\textquotedbl{} \cite{Hawking2006}
and many others \cite{Hodges1998,Aczel2000}; this version can be
forced to fail, and the same happens with D. R. Hofstadter's variation
\cite{Hofstadter1999} that subtracts 1 from every digit modulo 10.
On the other hand, R. Penrose's version \cite{Penrose1991} cannot
be forced to yield a number on the list, and the resulting antidiagonal
number is an unwritable rational. W. Dunham's approach \cite{Dunham1992}
uses random digits, which can neither be forced to fail (but note
that random numbers cannot be generated by a purely mathematical process---they
require an entropy source, which is a physical device).

\begin{table}[htbp]
\centering{}\iftoggle{epstablesarenicer}{\includegraphics[bb=72bp 612bp 515bp 708bp,clip,width=1\textwidth]{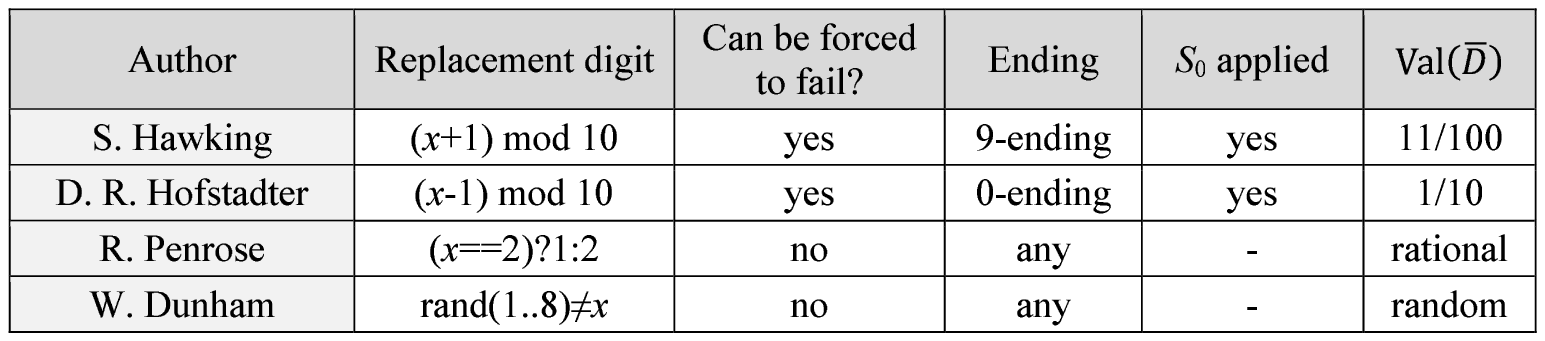}}{%
\begin{tabular}{|>{\raggedright}p{80pt}|>{\centering}p{65pt}|>{\centering}p{65pt}|c|>{\centering}p{33pt}|c|}
\hline 
\rowcolor{gray09}Author & Replacement\\
 digit & Can be forced\\
 to fail? & Ending & $S_{0}$\\
applied & $\mathrm{Val}\left(\overline{D}\right)$\tabularnewline
\hline 
S. Hawking & $(x+1)\,\mathrm{mod}\,10$ & yes & 9-ending & yes & 11/100 \tabularnewline
\hline 
D. R. Hofstadter & $(x-1)\,\mathrm{mod}\,10$ & yes & 0-ending & yes & 1/10 \tabularnewline
\hline 
R. Penrose & $(x=2)?\,1\,:\,2$ & no & any & - & rational \tabularnewline
\hline 
W. Dunham & $\mathrm{rand}(1..8)\ne x$ & no & any & - & random \tabularnewline
\hline 
\end{tabular}}\caption{Results of the application of the DM to $L_{\mathrm{DI}}$ in base
10 for some of the variants found in popular math vulgarization books.}
\label{tab:Table6} 
\end{table}

With regard to the replacement digit criterion, the base 2 case for
the DM is especially hopeless. Given that in base 2 there are only
two available symbols it is always possible to find a shuffling of
the list that forces a never-ending ...111... or ...000... trivial
termination for the antidiagonal number, which as we have seen converges
to numbers in the closure of $\mathbb{W}_{2}$ that are also representable
within $\mathbb{W}_{2}$. This situation can ultimately be forced
because the positional fractional notation is more efficient at producing
strings (numbers) than the DM is at discarding them. Considering that
$L_{\mathrm{DI}}$ only needs $\log_{2}\left(n\right)$ digits to
represent the number at line \emph{n} and the DM needs to swap digit
\emph{n} to clear the number at line \emph{n}, this implies the DM
will eventually find an area where there are no longer any significant
digits (meaning the trivial ...111... and ...000... infinite terminations
on the list). Consequently, the DM will---for infinitely many appropriate
shuffling functions of $L_{\mathrm{DI}}$, if needed---end in a constant
tail of 1s or 0s (depending on the flavor chosen for the DM), thus
converging into a number that can be represented by a finite number
of digits and will therefore be explicitly contained within $L_{\mathrm{DI}}$.

\section{Application of the DM to a generic list of reals\label{sec:Apply-Ldi-to-R}}

What if instead of applying the DM to $L\left([0,\thinspace1)_{\mathbb{W}_{2}}\right)$
we used $L\left([0,\thinspace1)\right)$ as an input? Would the DM
produce a noticeably different result? That could perhaps indicate
that the DM can tell apart lists of denumerable sets and list of (presumably)
non-denumerable sets--although it is very unlikely because the DM
black box (Figure \ref{fig:Figure1}) requires first to write down
the numbers using positional notation, and there are no more written
numbers than those in $L\left(\mathbb{W}\right)$.

In this section we will just demonstrate that some orderings of $L\left(\mathbb{R}\right)$
also force the DM to yield an antidiagonal number $\overline{D}$
that is already in $L\left(\mathbb{R}\right)$. The proof relies in
the previously mentioned fact that the DM cannot distinguish unwritable
real numbers (with infinitely many significant digits) from numbers
in $\mathbb{W}$, and therefore we can make sure our input list of
reals yields one of the numbers of the list by replicating, for instance,
the \textquotedbl{}skeleton\textquotedbl{} of the list $L_{\mathrm{DI}}^{''}$
shown in Table \ref{tab:Table5}---which we know forces the DM to
fail.

So we start with a list $L_{R}=L\left((0,\thinspace1]\right)$ which
we assume contains all real numbers in (0, 1{]} written using 1-ending
termination and proceed to construct a reordering of it, $L_{R}^{'}$,
in the following way:
\begin{enumerate}
\item Start by making $L_{R}^{'}=L_{R}$.
\item Find the number 3/4 (the real value of $0.101111\ldots_{2}$) in $L_{R}^{'}$
and move it to line 3, shifting the rest of the list as needed.
\item Starting from line 1, check if the \emph{k}\textsuperscript{th} binary
digit of $L_{R}^{'}\left(k\right)$ using 1-ending representation
matches that of $L_{\mathrm{DI}}^{''}\left(k\right)$. If it does,
do nothing; otherwise look ahead in $L_{R}^{'}$ until you find an
index \emph{m} such that $L_{R}^{'}\left(m\right)$ matches the \emph{k}\textsuperscript{th}
digit, then swap lines \emph{m} and \emph{k} in $L_{R}^{'}$.
\item Repeat step 3 until the whole list has been reordered.
\end{enumerate}
An example of application of this procedure is shown in Table \ref{tab:Table7}.

\begin{table}[htbp]
\centering{}\iftoggle{epstablesarenicer}{\includegraphics[bb=90bp 560bp 501bp 708bp,clip,width=1\textwidth]{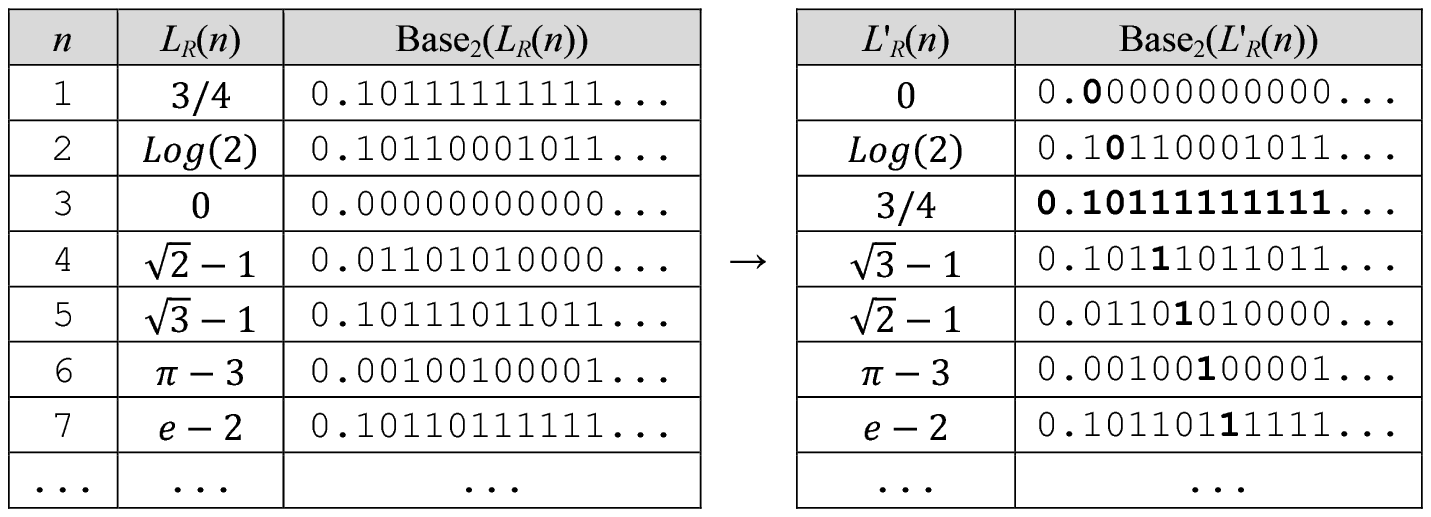}}{%
\begin{tabular}{|>{\centering}p{23pt}|>{\centering}p{35pt}|l|>{\raggedright}p{20pt}|>{\centering}p{35pt}|l|}
\cline{1-3} \cline{5-6} 
\rowcolor{gray09}\emph{n} & $L_{R}\left(n\right)$ & $\mathrm{Base_{2}}\left(L_{R}\left(n\right)\right)$ & \cellcolor{white} & $L_{R}^{'}\left(n\right)$ & $\mathrm{Base_{2}}\left(L_{R}^{'}\left(n\right)\right)$\tabularnewline
\cline{1-3} \cline{5-6} 
1 & $3/4$ & 0.10111111111... &  & $0$ & 0.\textbf{0}0000000000... \tabularnewline
\cline{1-3} \cline{5-6} 
2 & $\log\left(2\right)$ & 0.10110001011... &  & $\log\left(2\right)$ & 0.1\textbf{0}110001011... \tabularnewline
\cline{1-3} \cline{5-6} 
3 & $0$ & 0.00000000000... &  & $3/4$ & \textbf{0.10111111111...} \tabularnewline
\cline{1-3} \cline{5-6} 
4 & $\sqrt{2}\mathrm{-1}$ & 0.01101010000... & $\longrightarrow$ & $\sqrt{3}\mathrm{-1}$ & 0.101\textbf{1}1011011... \tabularnewline
\cline{1-3} \cline{5-6} 
5 & $\sqrt{3}\mathrm{-1}$ & 0.10111011011... &  & $\sqrt{2}\mathrm{-1}$ & 0.0110\textbf{1}010000... \tabularnewline
\cline{1-3} \cline{5-6} 
6 & $\pi-3$ & 0.00100100001... &  & $\pi-3$ & 0.00100\textbf{1}00001... \tabularnewline
\cline{1-3} \cline{5-6} 
7 & $e-2$ & 0.10110111111... &  & $e-2$ & 0.101101\textbf{1}1111... \tabularnewline
\cline{1-3} \cline{5-6} 
... & ... & ... &  & ... & ... \tabularnewline
\cline{1-3} \cline{5-6} 
\end{tabular}}\caption{A reordering $L_{R}^{'}$ of the list $L_{R}$ of the reals in {[}0,1)
so that its diagonal and third element match those of $L_{\mathrm{DI}}^{''}$.}
\label{tab:Table7} 
\end{table}

\begin{figure}[h]
\begin{centering}
\includegraphics[bb=74bp 440bp 528bp 705bp,clip,scale=0.5]{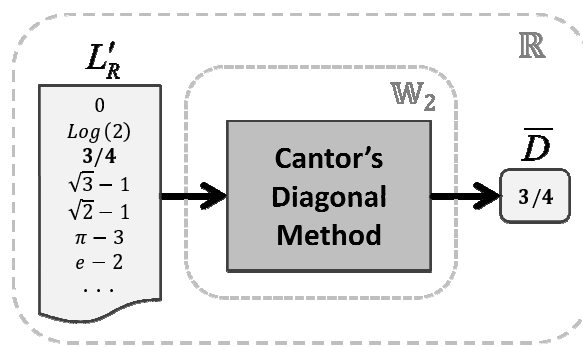}
\par\end{centering}

\caption{The list $L_{R}^{'}$ of the reals in {[}0,1), which has been specially
engineered to force Cantor's Diagonal Method to yield the third element
of the list when working on base 2. Note that the DM block must always
operate within a proper subset of $\mathbb{R}$, $\mathbb{W}_{b}$.}

\label{fig:Figure2}
\end{figure}

With this reordering of $L_{R}$ we are sure to obtain $L_{R}^{'}\left(3\right)$
when applying the binary DM, since $L_{R}^{'}$ and $L_{\mathrm{DI}}^{''}$
are indistinguishable to the DM. Figure \ref{fig:Figure2} summarizes
what we have achieved.

We can then conclude that there is at least one shuffling of $L_{R}$
that forces the DM to yield a number explicitly listed. Therefore
we have just proved that it is not true that for any given list of
the reals in (0,1{]} the DM is always able to produce a number not
on the list.

Also, we have shown that the representation of the input list plays
an essential role in the outcome of the procedure: in base 2 we can
force the DM to fail, but for other bases b>2\textemdash{}and some
careful selection of the replacement digits\textemdash{}it can still
perform as expected. We will see later that this is not a mere coincidence,
and that the cardinality of infinite sets is determined by the representation
we choose for modeling them rather than on their \textquotedbl{}size\textquotedbl{}.

But first, let us analyze the relevance of the DM to Cantor's Theorem,
which ultimately is what is assumed to guarantee the existence of
infinite cardinalities higher than $\aleph_{0}$.

\section{A counterexample to Cantor's Theorem for infinite sets\label{sec:Counter-Cantor-Theorem}}

The ultimate reason why the DM has no relation to the cardinality
of $\mathbb{R}$ is that Cantor's Theorem does not hold for infinite
sets. After all, the result of the DM alone would just be that one
single element is missing from the list, which by itself is not enough
to prove a difference in cardinality. In combination with Cantor's
Theorem it may however have some ground.

Cantor's Theorem establishes that it is not possible to put one set
in correspondence with its power-set. While this is certainly true
for finite sets, it's not at all obvious for infinite sets, since
in principle we have as many elements as we may need. A refutation
is therefore required, and a counterexample is probably the best way
to do it . To produce it we need to introduce some concepts first.

A selector or indicator function $\chi_{C}$ over a set \emph{A} permits
to choose a subset \emph{C} by taking the value 1 for the elements
that belong to \emph{C} and 0 for the rest: 
\begin{equation}
C\subset A
\end{equation}
\begin{equation}
\chi_{C}\thinspace:\thinspace A\to\left\{ 0,\thinspace1\right\} 
\end{equation}
\begin{equation}
\chi_{C}\left(x\right)=\left\{ \begin{array}{cc}
1 & if\thinspace x\in C\\
0 & if\thinspace x\notin C
\end{array}\right.
\end{equation}
If the original set \emph{A} is ordered in a list $L_{A}$ the indicator
function can take the shape of a binary string $w_{C}$ where the
position of its zeros and ones indicate the positions of those elements
in $L_{A}$ that belong to the subset \emph{C}: 
\begin{equation}
w_{C}\in\mathbb{W}_{2}^{I}
\end{equation}
\begin{equation}
w_{C}=x_{n}\ldots x_{1}x_{0}\sim\chi_{C}\Leftrightarrow x_{i}=\chi_{C}\left(L_{A}\left(i\right)\right)\thinspace\thinspace\forall i\in\left\{ 0,\ldots,n\right\} \label{eq:IndicatorEq}
\end{equation}
For example, we can choose the subset \{3, 22\} from the ordered set
\emph{L} = \{3, 42, 2, 22\} with the binary number \emph{c} = 1001.

Using this representation we can see that for a set \emph{A} of \emph{n}
elements there are $2^{n}$ possible subsets. The collection of all
those subsets is called the power set of \emph{A,} \emph{P}(\emph{A}).
Cantor's Theorem asserts that Card(\emph{A}) < Card(\emph{P}(\emph{A}))
for every non-empty set \emph{A}, which in practice means that it
is not possible to find a bijection between any set and its power-set,
including the case where \emph{A} is an infinite set. Applying the
theorem to the set of natural numbers $\mathbb{N}$, which has cardinality
$\aleph_{0}$, results in $\aleph_{0}<2^{\aleph_{0}}=\mathfrak{c}$,
where $\mathfrak{c}$ is the cardinality of the set of real numbers
$\mathbb{R}$. For Cantor this result implies that it is impossible
to create an ordered list that contains all real numbers, and thus
that $\mathbb{R}$ is a non-denumerable set.

The demonstration of Cantor's Theorem, however, uses the diagonal
argument, which we have already seen is not reliable in some situations.
We need to determine if application of the DM within the context of
Cantor's Theorem is justified or not.

We have already seen that $\mathbb{W}_{2}^{I}$ is isomorphic to $\mathbb{Z}_{2}^{+}$
and the equivalence shown in (\ref{eq:IndicatorEq}) seems to indicate
that the set of selectors for $\mathbb{Z}^{+}$ can be well-ordered.

Let us consider the set of non-negative integers and their corresponding
binary representation, shown in Table \ref{tab:Table8} below, where
$\mathrm{s}$ indicates the successor function. I'm writing the non-negative
numbers using the successor function to again remind the reader that
natural numbers, like the reals, don't have digits either---they are
abstract objects that, incidentally, are convenient to represent using
positional notation.

\begin{table}[htbp]
\centering{}\iftoggle{epstablesarenicer}{\includegraphics[bb=135bp 394bp 477bp 708bp,clip]{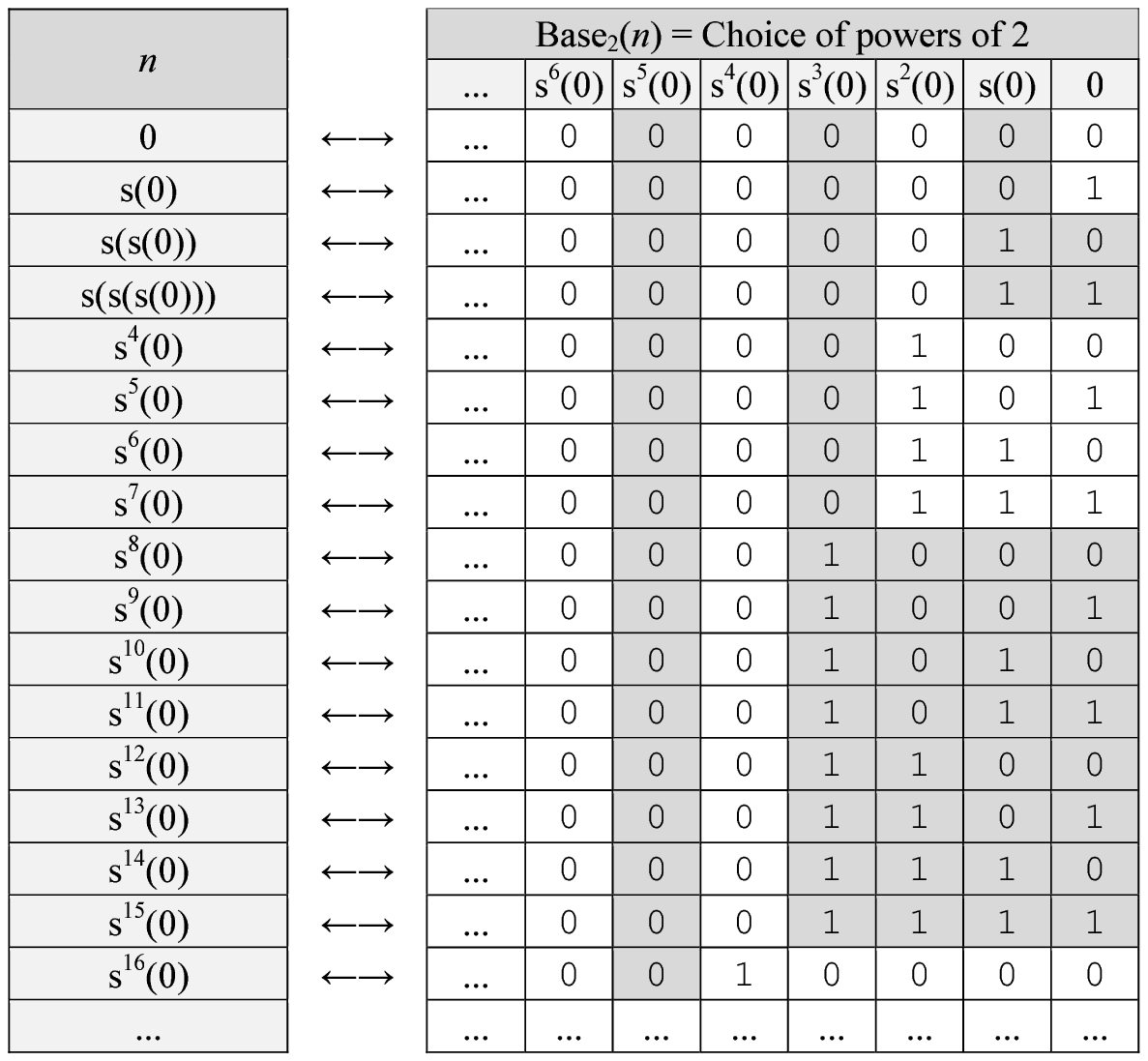}}{%
\begin{tabular}{|>{\centering}p{46pt}|>{\centering}p{25pt}|>{\centering}p{16pt}|>{\centering}p{23pt}|>{\centering}p{24pt}|>{\centering}p{24pt}|>{\centering}p{24pt}|>{\centering}p{24pt}|>{\centering}p{24pt}|>{\centering}p{24pt}|}
\multicolumn{1}{>{\centering}p{46pt}}{} & \multicolumn{1}{>{\centering}p{25pt}}{} & \multicolumn{8}{p{197pt}}{\centering$\mathrm{Base}_{2}(n)=Choice\, of\, powers\, of\,2$}\tabularnewline
\cline{1-1} \cline{3-10} 
\rowcolor{gray09}n & \cellcolor{white} & ... & $\mathrm{s}{}^{\mathrm{6}}\left(0\right)$ & $\mathrm{s}{}^{5}\left(0\right)$ & $\mathrm{s}{}^{4}\left(0\right)$ & $\mathrm{s}{}^{\mathrm{3}}\left(0\right)$ & $\mathrm{s}{}^{\mathrm{2}}\left(0\right)$ & $\mathrm{s}\left(0\right)$ & 0 \tabularnewline
\cline{1-1} \cline{3-10} 
0 & $\longleftrightarrow$ & ... & 0 & \cellcolor{gray097}0 & 0 & \cellcolor{gray097}0 & 0 & \cellcolor{gray097}0 & 0 \tabularnewline
\cline{1-1} \cline{3-10} 
$\mathrm{s}\left(0\right)$ & $\longleftrightarrow$ & ... & 0 & \cellcolor{gray097}0 & 0 & \cellcolor{gray097}0 & 0 & \cellcolor{gray097}0 & 1 \tabularnewline
\cline{1-1} \cline{3-10} 
$\mathrm{s}\left(\mathrm{s}\left(0\right)\right)$ & $\longleftrightarrow$ & ... & 0 & \cellcolor{gray097}0 & 0 & \cellcolor{gray097}0 & 0 & \cellcolor{gray097}1 & \cellcolor{gray097}0 \tabularnewline
\cline{1-1} \cline{3-10} 
$\mathrm{s}\left(\mathrm{s}\left(\mathrm{s}\left(0\right)\right)\right)$ & $\longleftrightarrow$ & ... & 0 & \cellcolor{gray097}0 & 0 & \cellcolor{gray097}0 & 0 & \cellcolor{gray097}1 & \cellcolor{gray097}1 \tabularnewline
\cline{1-1} \cline{3-10} 
$\mathrm{s}{}^{4}\left(0\right)$ & $\longleftrightarrow$ & ... & 0 & \cellcolor{gray097}0 & 0 & \cellcolor{gray097}0 & 1 & 0 & 0 \tabularnewline
\cline{1-1} \cline{3-10} 
$\mathrm{s}{}^{5}\left(0\right)$ & $\longleftrightarrow$ & ... & 0 & \cellcolor{gray097}0 & 0 & \cellcolor{gray097}0 & 1 & 0 & 1 \tabularnewline
\cline{1-1} \cline{3-10} 
$\mathrm{s}{}^{\mathrm{6}}\left(0\right)$ & $\longleftrightarrow$ & ... & 0 & \cellcolor{gray097}0 & 0 & \cellcolor{gray097}0 & 1 & 1 & 0 \tabularnewline
\cline{1-1} \cline{3-10} 
$\mathrm{s}{}^{\mathrm{7}}\left(0\right)$ & $\longleftrightarrow$ & ... & 0 & \cellcolor{gray097}0 & 0 & \cellcolor{gray097}0 & 1 & 1 & 1 \tabularnewline
\cline{1-1} \cline{3-10} 
$\mathrm{s}{}^{\mathrm{8}}\left(0\right)$ & $\longleftrightarrow$ & ... & 0 & \cellcolor{gray097}0 & 0 & \cellcolor{gray097}1 & \cellcolor{gray097}0 & \cellcolor{gray097}0 & \cellcolor{gray097}0 \tabularnewline
\cline{1-1} \cline{3-10} 
$\mathrm{s}{}^{\mathrm{9}}\left(0\right)$ & $\longleftrightarrow$ & ... & 0 & \cellcolor{gray097}0 & 0 & \cellcolor{gray097}1 & \cellcolor{gray097}0 & \cellcolor{gray097}0 & \cellcolor{gray097}1 \tabularnewline
\cline{1-1} \cline{3-10} 
$\mathrm{s}{}^{\mathrm{10}}\left(0\right)$ & $\longleftrightarrow$ & ... & 0 & \cellcolor{gray097}0 & 0 & \cellcolor{gray097}1 & \cellcolor{gray097}0 & \cellcolor{gray097}1 & \cellcolor{gray097}0 \tabularnewline
\cline{1-1} \cline{3-10} 
$\mathrm{s}{}^{\mathrm{11}}\left(0\right)$ & $\longleftrightarrow$ & ... & 0 & \cellcolor{gray097}0 & 0 & \cellcolor{gray097}1 & \cellcolor{gray097}0 & \cellcolor{gray097}1 & \cellcolor{gray097}1 \tabularnewline
\cline{1-1} \cline{3-10} 
$\mathrm{s}{}^{\mathrm{12}}\left(0\right)$ & $\longleftrightarrow$ & ... & 0 & \cellcolor{gray097}0 & 0 & \cellcolor{gray097}1 & \cellcolor{gray097}1 & \cellcolor{gray097}0 & \cellcolor{gray097}0 \tabularnewline
\cline{1-1} \cline{3-10} 
$\mathrm{s}{}^{\mathrm{13}}\left(0\right)$ & $\longleftrightarrow$ & ... & 0 & \cellcolor{gray097}0 & 0 & \cellcolor{gray097}1 & \cellcolor{gray097}1 & \cellcolor{gray097}0 & \cellcolor{gray097}1 \tabularnewline
\cline{1-1} \cline{3-10} 
$\mathrm{s}{}^{\mathrm{14}}\left(0\right)$ & $\longleftrightarrow$ & ... & 0 & \cellcolor{gray097}0 & 0 & \cellcolor{gray097}1 & \cellcolor{gray097}1 & \cellcolor{gray097}1 & \cellcolor{gray097}0 \tabularnewline
\cline{1-1} \cline{3-10} 
$\mathrm{s}{}^{\mathrm{15}}\left(0\right)$ & $\longleftrightarrow$ & ... & 0 & \cellcolor{gray097}0 & 0 & \cellcolor{gray097}1 & \cellcolor{gray097}1 & \cellcolor{gray097}1 & \cellcolor{gray097}1 \tabularnewline
\cline{1-1} \cline{3-10} 
$\mathrm{s}{}^{16}\left(0\right)$ & $\longleftrightarrow$ & ... & 0 & \cellcolor{gray097}0 & 1 & 0 & 0 & 0 & 0 \tabularnewline
\cline{1-1} \cline{3-10} 
... &  & ... & ... & ... & ... & ... & ... & ... & ... \tabularnewline
\cline{1-1} \cline{3-10} 
\end{tabular}}\caption{A bijection between the set of non-negative integers \emph{n} in $\mathbb{Z}^{+}$
and their binary representation, which actually matches the set of
choice functions of the integer powers of 2 whose sum is \emph{n}.
The shaded areas delineate the $2^{k}$ blocks that contain all possible
choices over the first \emph{k} integers; when all $\aleph_{0}$ integers
are used the block corresponds to the whole list which must therefore
contain $2^{\aleph_{0}}$ elements.}
\label{tab:Table8} 
\end{table}

Table \ref{tab:Table8} establishes a bijection between any number
in $\mathbb{Z}^{+}$ and its binary representation string in $\mathbb{W}_{2}^{I}=\mathbb{Z}_{2}^{+}$.
Now, the binary representation is just a choice function of the different
powers of 2 which are themselves the numbers in $\mathbb{Z}^{+}$;
every digit 0/1 indicates if the number signaling the power is included
in the choice or not. This basically means that the power-set of $\mathbb{Z}^{+}$
is equipollent to $\mathbb{Z}^{+}$, and therefore their cardinalities
have to be equivalent. Since we can use as many as $\aleph_{0}$ digits,
the set of binary indicators has cardinality $2^{\aleph_{0}}$; but
we know the set of binary integer numbers also has cardinality $\aleph_{0}$
because it is $\mathbb{Z}^{+}$ itself. Therefore we must conclude
that $2^{\aleph_{0}}=\aleph_{0}$, implying at the same time that
Cantor's Theorem can't possibly hold for infinite sets and that the
cardinality of the continuum $\mathfrak{c}$ is likely to be equivalent
to that of the natural numbers $\aleph_{0}$.

Note that this actually shows that the representation we choose for
$\mathbb{Z}^{+}$ is what ultimately determines its cardinality type:
using the successor representation we get cardinality $\aleph_{0}$,
while using binary positional notation we get $2^{\aleph_{0}}$; and
using a factorial numbering system \cite{Laisant1888} we would get
cardinality $\aleph_{0}!$. The set $\mathbb{Z}^{+}$ has still the
same ``size'' of $\infty$ in all three cases. A similar situation
happens with the set of rational numbers $\mathbb{Q}$: using the
positional fractional numbering representation yields a cardinality
of $2^{\aleph_{0}}$ for $\mathbb{Q}$ while the p/q representation
results in a cardinality of $\mathrm{Card}\left(\mathbb{Z}\times\mathbb{N}\right)=\aleph_{0}$;
$\mathbb{Q}$ does not suffer any variation in size because of the
change of representation. The relevance of the representation in the
determination of the cardinality of a set is covered in detail on
the next section.

Note as well that once we have provided such an obvious bijection
(basically from $\mathbb{Z}^{+}$ to itself) between an ordered infinite
set and its power-set, we don't even need to discuss what may be the
results of the application of the DM to the selector list in Table
\ref{tab:Table8}. In any case, the application yields the antidiagonal
object $\bar{\Delta}=\ldots1111_{2}$ which is usually regarded as
the selector for the whole set $\mathbb{Z}^{+}$. However, we can
see this not to be the case, since it selects more than integer numbers:

\begin{equation}
\bar{\Delta}=\ldots1111_{2}=\lim_{n\rightarrow\infty}\sum\limits _{k=0}^{n}{1\cdot2^{k}}=\sum\limits _{k=0}^{\infty}{1\cdot2^{k}}=1\cdot2^{\infty}+\cdots+1\cdot2^{0}
\end{equation}

And we can see that this potential selector is also choosing infinity,
which is not a member of the set of integers. Therefore, whatever
the object $\bar{\Delta}$ happens to be, it is not the selector for
all the integers. We need to remember at this point that objects in
the closure of a set may not belong to the set (the limits of converging
rational sequences may be irrational for instance), which is precisely
the case here.

In fact, we can see that there is no selector for the whole $\mathbb{Z}^{+}$
by way of a contradiction: assume we have a selector for all $z\in\mathbb{Z}^{+}$;
therefore, since we have selected all of them and because they are
ordered, there must be a last element \emph{w} that we have selected;
however, since $w\in\mathbb{Z}^{+}$ so does $\mathrm{S}\left(w\right)=w+1$;
this means our selector hasn't selected all numbers. It is impossible
to select all numbers in $\mathbb{Z}^{+}$ and so the selector for
the whole of $\mathbb{Z}^{+}$ simply does not exist.

This points to an interesting characteristic of infinite sets: they
can be defined or described, but not completely realized. It is easy
to determine if an object \emph{x} belongs to $\mathbb{Z}^{+}$ just
by checking if it satisfies the properties of the elements in $\mathbb{Z}^{+}$,
but enumerating or listing every object that may satisfy them is not
possible. This would be a trivial remark were it not for completed
infinity being an essential part of transfinite number theory.

We could counter-argue that $\mathbb{W}_{2}^{I}$ does not cover all
possible selections that may exist, but that is not sufficient to
invalidate the fact that $\aleph_{0}=2^{\aleph_{0}}$. In Table \ref{tab:Table8}
we have an ordered set with $\aleph_{0}$ elements paired to another
ordered set that has $2^{\aleph_{0}}$ elements, and this pairing
is surely a bijection because we are matching $\mathbb{Z}^{+}$ to
itself (actually we are matching one representation of $\mathbb{Z}^{+}$
to another representation of $\mathbb{Z}^{+}$, but we will show in
section \ref{sec:Numbering-models} that the bijection is preserved
when we demonstrate that both the successor and positional representations
form proper models for $\mathbb{Z}^{+}$). This is similar to pairing
off the set of rationals to the even numbers to prove the rationals
have cardinality $\aleph_{0}$: we may have left some numbers out
somewhere, but both sets have preserved their cardinalities in the
pairing.

\subsection{The Applicative Numbering Model\label{sub:The-Applicative-Numbering}}

To further ratify that differences in cardinality are not the main
reason for a set being denumerable or not, we will show another representation
for the integers that also has cardinality $\mathfrak{c}$: consider
the limit case where a positional notation system with $\aleph_{0}$
available symbols (number base $b=\infty$) is used to represent $\mathbb{N}$.
In this situation we can construct up to $\aleph_{0}^{\aleph_{0}}$
different strings of finite length. We will now show that such a set
can be well-ordered.

It is evident that well-ordering the set $\mathbb{N}_{\aleph_{0}}$
of $\aleph_{0}^{\aleph_{0}}$ elements is not possible using the regular
lexicographic ordering because we would exhaust the set of indices
$\mathbb{Z}^{+}$ before reaching two digits' strings; the lexicographic
order is therefore not a suitable mapping function for a proper model
of such a set (proper numbering models are covered in full on the
next section).

In order to tackle the problem we need to ensure that every isolated
symbol will appear at a finite position and also that all possible
combinations of any number of symbols will be present in the list.
A potential approach is to group the strings in blocks that use up
to \emph{n} digits using up to \emph{n} symbols. This results in a
recursive structure that introduces new symbols only when required
while keeping the string length growing steadily.

The block structure for this Applicative Numbering model can be generated
by reordering the standard lexicographic representation so that every
string that involves either the last \emph{n}\textsuperscript{th}
symbol or requires \emph{n} digits is moved to the end of the block;
the initial part of the block is the applicative numbering for up
to \emph{n}-1 digits using up to \emph{n}-1 symbols, thus ensuring
the recursive structure needed. Table \ref{tab:Table9} below shows
the applicative block structure for the strings of up to 3 digits
using up to 3 symbols.

\begin{table}[htbp]
\centering{}\iftoggle{epstablesarenicer}{\includegraphics[bb=119bp 364bp 335bp 708bp,clip]{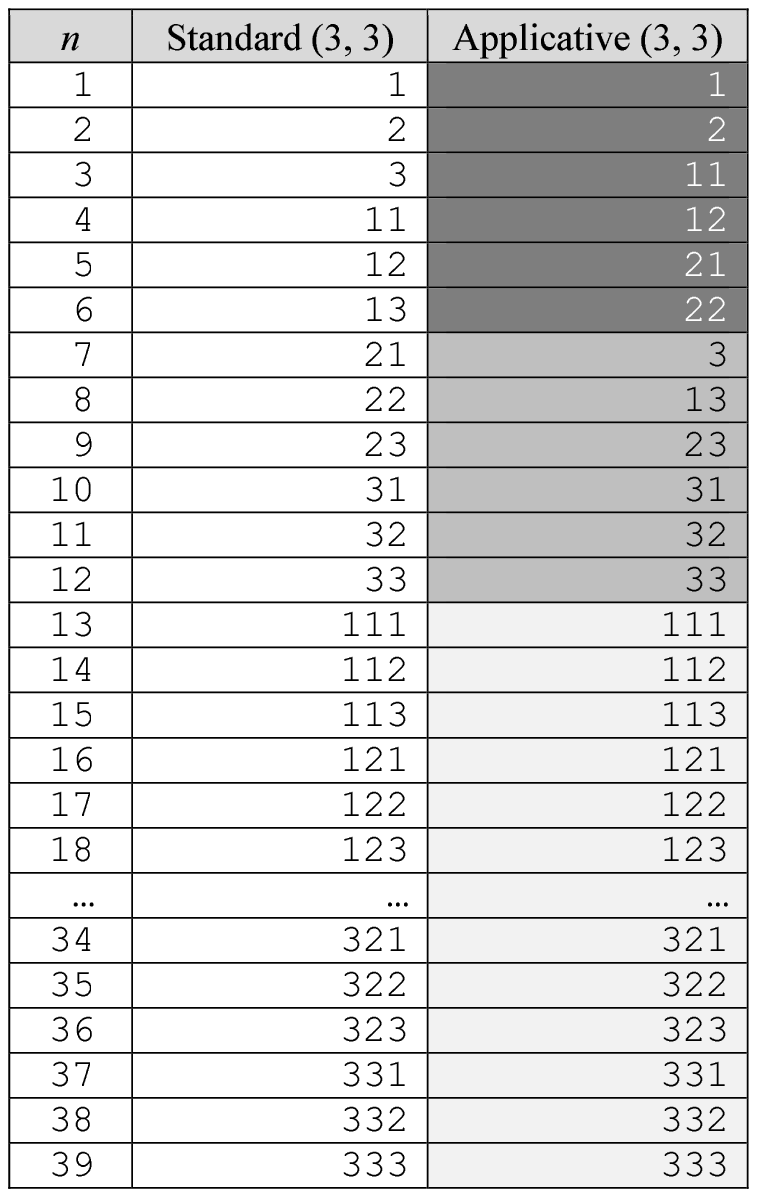}}{%
\begin{tabular}{|>{\centering}p{35pt}|r|r|}
\hline 
\rowcolor{gray09}n & Standard (3, 3) & Applicative (3, 3) \tabularnewline
\hline 
1 & 1 & \textcolor{white}{\cellcolor{gray033}1 }\tabularnewline
\hline 
2 & 2 & \textcolor{white}{\cellcolor{gray033}2 }\tabularnewline
\hline 
3 & 3 & \textcolor{white}{\cellcolor{gray033}11 }\tabularnewline
\hline 
4 & 11 & \textcolor{white}{\cellcolor{gray033}12 }\tabularnewline
\hline 
5 & 12 & \textcolor{white}{\cellcolor{gray033}21 }\tabularnewline
\hline 
6 & 13 & \textcolor{white}{\cellcolor{gray033}22 }\tabularnewline
\hline 
7 & 21 & \cellcolor{gray09}3 \tabularnewline
\hline 
8 & 22 & \cellcolor{gray09}13 \tabularnewline
\hline 
9 & 23 & \cellcolor{gray09}23 \tabularnewline
\hline 
10 & 31 & \cellcolor{gray09}31 \tabularnewline
\hline 
11 & 32 & \cellcolor{gray09}32\tabularnewline
\hline 
12 & 33 & \cellcolor{gray09}33 \tabularnewline
\hline 
13 & 111 & \cellcolor{gray097}111 \tabularnewline
\hline 
14 & 112 & \cellcolor{gray097}112 \tabularnewline
\hline 
15 & 113 & \cellcolor{gray097}113 \tabularnewline
\hline 
16 & 121 & \cellcolor{gray097}121 \tabularnewline
\hline 
17 & 122 & \cellcolor{gray097}122 \tabularnewline
\hline 
18 & 123 & \cellcolor{gray097}123 \tabularnewline
\hline 
$\ldots$  & $\ldots$ & \cellcolor{gray097}$\ldots$\tabularnewline
\hline 
34 & 321 & \cellcolor{gray097}321 \tabularnewline
\hline 
35 & 322 & \cellcolor{gray097}322 \tabularnewline
\hline 
36 & 323 & \cellcolor{gray097}323 \tabularnewline
\hline 
37 & 331 & \cellcolor{gray097}331 \tabularnewline
\hline 
38 & 332 & \cellcolor{gray097}332 \tabularnewline
\hline 
39 & 333 & \cellcolor{gray097}333 \tabularnewline
\hline 
\end{tabular}}\caption{Standard lexicographic ordering and applicative ordering of all the
strings with up to 3 digits using 3 symbols.}
\label{tab:Table9} 
\end{table}

The shading in the right column of Table \ref{tab:Table9} shows that
the applicative block for (\emph{n}, \emph{n}) is composed of three
main parts: the initial one is the applicative block for (\emph{n}-1,
\emph{n}-1) and the last one---which uses exactly \emph{n} digits---matches
that of the standard lexicographic ordering; the middle part of the
(\emph{n}, \emph{n}) block contains all strings that require less
than \emph{n} digits but contain the \emph{n}\textsuperscript{th}
symbol at least once, and should therefore be after the (\emph{n}-1,
\emph{n}-1) part. Note as well that the (\emph{n}-1, \emph{n}-1) block
also follows the same recursive pattern.

Working in this way we can see that the addition of another symbol
will not interfere with the previously listed elements, and also that
strings of any length are eventually reached.

The cardinality of the applicative set of up to \emph{d} digits using
\emph{n} symbols is the number of applications from the set of symbols
to the set of digits:

\begin{equation}
A_{n}^{[1,\thinspace d]}=n+n^{2}+\cdots+n^{d-1}+n^{d}=\frac{n^{d+1}-1}{n-1}-1
\end{equation}
And the limit when using up to $\aleph_{0}$ symbols and $\aleph_{0}$
digits is: 
\begin{equation}
\lim_{n\rightarrow\aleph_{0}}A_{n}^{[1,n]}=\frac{\aleph_{0}^{\aleph_{0}+1}-1}{\aleph_{0}-1}-1=\frac{\aleph_{0}^{\aleph_{0}}}{\aleph_{0}}-1=\aleph_{0}^{\aleph_{0}-1}-1=\aleph_{0}^{\aleph_{0}}
\end{equation}
And we can see that with the applicative numbering we can well-order
sets of cardinality $\aleph_{0}^{\aleph_{0}}$. Since we know that
$\aleph_{0}^{\aleph_{0}}=2^{\aleph_{0}}=\mathfrak{c}$ \cite{Sierpinski1965},
we have therefore shown again that a set of cardinality $\mathfrak{c}$
can be well-ordered.

\section{General considerations for well-ordering a set\label{sec:Numbering-models}}

To order a set \emph{S} we need to find a bijection between it and
the set of non-negative integers. What is often overlooked is that
the bijection is not between \emph{S} and $\mathbb{Z}^{+}$, but between
a representation of \emph{S}, \emph{R}(\emph{S}), and a representation
of $\mathbb{Z}^{+}$, $R\mathbb{\left(Z^{\mathrm{+}}\right)}$. This
translates into restrictions on the representations used, which must
be able to express any element of the sets being modeled.

In general, we normally define a set using a group of properties or
axioms. This declarative definition of the set does not in principle
have to produce any elements of the set (except perhaps the identity
and null elements). For this task we normally rely on a model \emph{M}(\emph{S})
of the set \emph{S} that constructs representations of the elements
of \emph{S}. For instance, an abstract set such as $\mathbb{R}$ can
be defined by declaring the properties that its members must satisfy;
then we use the standard model that uses positional fractional notation
to generate real numbers in the form of strings of digits. We now
need to determine if this standard notation has any shortcomings for
representing the set of real numbers.

In general we can assume the model \emph{M} of a set \emph{S} to be
composed of a representation \emph{R}, a constructor function Rep
and a value function Val. The representation \emph{R}(\emph{S}) describes
members of \emph{S} using strings formed over an alphabet $\mathrm{\Sigma}$
of symbols or characters, and is therefore a subset of the Kleene
closure $\mathrm{\Sigma}^{\mathrm{\ast}}$. The constructor function
Rep generates the corresponding representation in \emph{R} of an element
in \emph{S}, while the value function Val translates those representations
back into \emph{S}.

\begin{equation}
M\left(S\right)=\left\{ R\left(S\right);\thinspace\mathrm{Rep};\thinspace\mathrm{Val};\thinspace\mathrm{\Sigma}\right\} 
\end{equation}
\[
R\left(S\right)\subset\mathrm{\Sigma}^{\mathrm{\ast}}
\]
\[
\mathrm{Rep}:S\to R(S)
\]
\[
\mathrm{Val}:R(S)\to S
\]
The model must obviously preserve the values of the elements of \emph{S},
and for this the Rep and Val functions must satisfy $\mathrm{Val}(\mathrm{Rep}(x))=x$
for all $x\in S$. However, it is not necessary for $\mathrm{Rep}(x)$
to be unique; that is, an element \emph{x} may have more than one
representation; for instance, consider 2/3 and 4/6, which represent
the same rational number 0.66666...

We say the model \emph{M}(\emph{S}) is valid, $M\left(S\right)\thinspace\overset{_{\mathrm{v}}}{\models}\thinspace S$,
if every element \emph{x} in \emph{S} has a representation Rep(\emph{x})
in the closure of \emph{R}(\emph{S}) that preserves its value.

\begin{equation}
M\left(S\right)\thinspace\overset{_{\mathrm{v}}}{\models}\thinspace S\thinspace\Longleftrightarrow\thinspace\forall x\in S\thinspace\exists r=\mathrm{Rep}\left(x\right)\in\mathrm{Cl}\left(R(S)\right)\thinspace\vert\thinspace\mathrm{Val}\left(r\right)=x
\end{equation}
We say the model \emph{M}(\emph{S}) is exhaustive, $M\left(S\right)\thinspace\overset{_{\mathrm{e}}}{\models}\thinspace S$,
if the values of the elements in the boundary of the representation
$\partial R\left(S\right)=\mathrm{Cl}\left(R(S)\right)\backslash R(S)$
don't belong to \emph{S} or, for those that do belong to \emph{S}
there is already an equivalent representation within \emph{R}(\emph{S}).

\begin{equation}
M\left(S\right)\thinspace\overset{_{\mathrm{e}}}{\models}\thinspace S\thinspace\Longleftrightarrow\thinspace\forall r\in\partial R\left(S\right):\thinspace\mathrm{Val}\left(r\right)\notin S\thinspace\vee\thinspace\exists r^{'}\in R\left(S\right)\thinspace\vert\thinspace\mathrm{Val}\left(r^{'}\right)=\mathrm{Val}\left(r\right)\in S
\end{equation}
Finally, we say that \emph{M} is a proper model of S, $M\left(S\right)\thinspace\overset{_{\mathrm{p}}}{\models}\thinspace S$,
if \emph{M}(\emph{S}) is both valid and exhaustive.

\begin{equation}
M\left(S\right)\thinspace\overset{_{\mathrm{p}}}{\models}\thinspace S\thinspace\Longleftrightarrow\thinspace\left(M\left(S\right)\thinspace\overset{_{\mathrm{v}}}{\models}\thinspace S\right)\thinspace\wedge\thinspace\left(M\left(S\right)\thinspace\overset{_{\mathrm{e}}}{\models}\thinspace S\right)
\end{equation}

Once we find a proper model for a set \emph{S}, it is straightforward
to find a well ordering for \emph{S}, because \emph{R}(\emph{S}) is
always a set of strings that can be lexicographically (or applicatively)
sorted into a list $L_{\mathrm{R(S)}}=L\left(R\left(S\right)\right)$.
The model being exhaustive means $L_{\mathrm{R(S)}}$ can then be
paired with $\mathbb{Z}^{+}$ without missing any elements in \emph{S},
because when $z\rightarrow\infty$ the objects $\mathrm{Val}\left(L_{\mathrm{R(S)}}\left(z\right)\right)$
either no longer belong to \emph{S} or have a previous equivalent
representation for some finite \emph{z'}. Any duplicated elements
in \emph{R}(\emph{S}) can later be trimmed into a simplified \emph{R'}(\emph{S})
by selecting the representation that comes first in $L_{\mathrm{R(S)}}$,
leaving a bijection between $\mathbb{Z}^{+}$ and \emph{R'}(\emph{S}),
and therefore between $\mathbb{Z}^{+}$ and \emph{S}.

We have seen that having a proper model for a set \emph{S} means that
the set can be well-ordered. Now we will see if the converse is true.
If a set \emph{S} accepts a well-ordering through a function $F:\mathbb{Z}^{+}\leftrightarrow S$
then we can use $\mathbb{Z}_{b}^{+}$ for any base \emph{b} as the
representation of \emph{S}, with $\mathrm{Rep}\left(x\right)=F^{-1}\left(x\right)$
as the constructor function and $\mathrm{Val}\left(x\right)=F\left(x\right)$
as the value function: $M\left(S\right)=\left\{ \mathbb{Z}_{b}^{+};\, F^{-1};\, F;\,\left\{ 0,\ldots,b-1\right\} \right\} $.
Now, since \emph{F} is a bijection there are no elements in \emph{S}
without a representation, which means \emph{M}(\emph{S}) is valid.
Also any elements in the closure $\lim_{z\rightarrow\infty}F\left(z\right)$
must either not belong to \emph{S} or be already represented for some
finite \emph{z'}, because otherwise \emph{F} wouldn't be a bijection;
this means \emph{M}(\emph{S}) is exhaustive.

We can thus justify the following theorem:

\begin{equation}
\exists L\left(S\right)\Leftrightarrow\exists M\thinspace\vert\thinspace M\left(S\right)\thinspace\overset{_{\mathrm{p}}}{\models}\thinspace S
\end{equation}
This theorem can be used as a constructive framework for Zermelo's
Well-Ordering Theorem, reducing the proposition that any set accepts
a well-ordering to an equivalent one that asserts that any set accepts
a proper model.

Also, once we have a proper model \emph{M} for a set \emph{S} we can
use \emph{M}(\emph{S}) to represent \emph{S} for any set-theoretic
operations---e.g. finding bijections between \emph{S} and other sets.

We can see now, for instance, that the model of the non-negative numbers
$\mathbb{Z}^{+}$ using concatenations of the successor function is
valid because every element in $\mathbb{Z}^{+}$ can be represented
in such a way. It is also exhaustive because the only element in the
limiting boundary $\lim_{n\rightarrow\infty}\mathrm{s}^{n}\left(0\right)=\infty$
does not belong to $\mathbb{Z}^{+}$. It can therefore be used to
well-order $\mathbb{Z}^{+}$.

\subsection{Proving $\mathbb{W}_{2}^{I}$ is a proper model for the set of selectors
of $\mathbb{Z}^{+}$\label{sub:Proving-WI2-proper}}

In order to justify the conclusions found in section \ref{sec:Counter-Cantor-Theorem}
regarding Cantor's Theorem it is important to prove that $\mathbb{W}_{2}^{I}=\mathbb{Z}_{2}^{+}$
is a proper model for the set of indicators of $\mathbb{Z}^{+}$ which
as we have seen is equivalent to the power-set of $\mathbb{Z}^{+}$,
$P\left(\mathbb{Z}^{+}\right)$. For this we need to show that it
has the properties of being valid and exhaustive.

We can see the model is valid because any collection \emph{C} of elements
from $\mathbb{Z}^{+}$ would result in a binary selector with 1s in
finite positions only (because the elements of $\mathbb{Z}^{+}$ are
all finite numbers), and would therefore be equivalent to a finite
non-negative binary integer (proof: it is bounded above by $2^{\left\lfloor \mathrm{log_{2}}\left(\mathrm{max}\left(C\right)\right)\right\rfloor +1}$
which is finite). This means the selector for any subset $C\subset\mathbb{Z}^{+}$
will be equivalent to a finite binary number in $\mathbb{Z}_{2}^{+}$,
and therefore $\mathbb{Z}_{2}^{+}$ is a valid model for the set of
indicators of $\mathbb{Z}^{+}$.

To see the model is exhaustive we can see that elements in the closure
of $\mathbb{Z}_{2}^{+}$ are only valid selectors if they start with
an infinite head of zeros; this makes them belong to the interior
of the closure. Limiting indicators in the boundary of the closure
are not valid indicators because they include 1s in non-finite positions,
thus selecting infinity which is not an element of $\mathbb{Z}^{+}$.
As a corollary, we can also see here that it is not possible to select
an infinity of only integers because the indicator number for \emph{C}
is bounded below by $2^{\mathrm{Card}\left(C\right)}-1$ which means
that any selection of an infinity of integers would be selecting infinity
as well.

Thus, $\mathbb{Z}_{2}^{+}$ is a proper model for the set of indicators
of $\mathbb{Z}^{+}$ and therefore a proper model for the power-set
$P\left(\mathbb{Z}^{+}\right)$. This justifies the conclusions derived
in section \ref{sec:Counter-Cantor-Theorem} on the applicability
of Cantor's Theorem to denumerable infinite sets.

\subsection{The standard positional numbering model\label{sub:Standard-numbering-model}}

We can also see that the standard positional fractional notation model
$M_{std}$ shown below is valid for $\mathbb{N}$, $\mathbb{Z}^{+}$,
$\mathbb{Q}^{+}$ and $\mathbb{R}^{+}$, because for any element of
those sets there is a positional fractional string of digits in the
standard model. In the case of $\mathbb{N}$ and $\mathbb{Z}^{+}$
the model is exhaustive, because the only element in the boundary
of the set, $\infty$, does not belong to them.

\begin{equation}
M_{std}\left(S\right)=\left\{ \begin{array}{l}
R_{std}\left(S\right)\equiv\mathrm{Cl}\left(\mathbb{W}_{b}\right);\thinspace\\
\mathrm{Rep_{\mathit{std}}}\left(x\in S\right)\equiv\left\{ w_{i}=\left\lfloor x\cdot b^{-i}\right\rfloor \thinspace mod\thinspace b\right\} _{i=-\infty}^{i=\infty};\thinspace\\
\mathrm{Val}_{std}\left(w=\left\{ w_{i}\right\} \in R_{std}\right)\equiv\sum\nolimits _{i=-\infty}^{i=\infty}{w_{i}\cdot b^{i}};\thinspace\\
\mathrm{\Sigma=}\left\{ \mathrm{0,\thinspace\ldots,\thinspace}b-1\mathrm{,\thinspace"."}\right\} 
\end{array}\right\} 
\end{equation}
In the case of $\mathbb{Q}$ and $\mathbb{R}$ the standard model
is not exhaustive, because there are elements that require an infinity
of digits in their representation and are therefore in the limiting
boundary of the representation $\partial R_{std}$. For the rationals
the situation can be easily corrected by augmenting (or replacing)
the standard representation with a new one $R_{pq}\left(\mathbb{Q}\right)$
which includes strings of the form \emph{p}/\emph{q}, where \emph{p}
and \emph{q} use the standard model for the integer and natural numbers
respectively. With this change we now have finite representations
equivalent to those infinite ones in the boundary that have a periodic
tail of fractional digits (e.g. 2/3 or 4/6 for 0.6666... ). The remaining
elements in the boundary are irrational numbers, that don't belong
to $\mathbb{Q}$, so this extended model is now exhaustive. With the
\emph{p}/\emph{q} representation we can then well-order the rationals
as Cantor did.

For $\mathbb{R}$ the situation is more complicated. Some real numbers
in the boundary of the representation $\partial R_{std}\left(\mathbb{R}\right)$
can be accounted for by introducing functions in the representation
such as square roots or logarithms, but there are transcendental numbers
that seem to escape simple algebraic reductions. We can therefore
state the important result:

\begin{equation}
\mathrm{\neg}\left(M_{std}\left(\mathbb{R}\right)\thinspace\overset{_{\mathrm{e}}}{\models}\thinspace\mathbb{R}\right)\Rightarrow\nexists L\left(R_{std}\left(\mathbb{R}\right)\right)
\end{equation}
That is, the standard positional representation of the real numbers
can't be well-ordered. This is what Cantor actually demonstrated.
From that, however, he incorrectly inferred that the set of reals
cannot be well-ordered. The inference would be justified if the standard
model was the only possible model for the reals. In reality however:

\begin{equation}
\nexists L\left(R_{std}\left(\mathbb{R}\right)\right)\nRightarrow\nexists L\left(\mathbb{R}\right)
\end{equation}
because there might be other alternative representations available
for $\mathbb{R}$, as we will see in the next section.

We can now understand why Cantor's diagonal argument seemed to imply
that the cardinality of $\mathbb{R}$ was greater than that of $\mathbb{N}$,
even if both sets are equally infinite: it was just because standard
model is not exhaustive for $\mathbb{R}$; therefore, any bijection
\emph{F} between $\mathbb{Z}^{+}$ and the standard representation
of $\mathbb{R}$ will always leave out some elements in the boundary
of the representation. Cantor's diagonal process---as well as his
1874 nested interval formulation---are actually ways to fetch one
of those elements left out by the chosen \emph{F}.

We can now also understand that infinite sets may have several valid
and exhaustive representations with different cardinalities: for instance,
$\mathbb{Z}^{+}$ has cardinality $\aleph_{0}$ when represented using
the successor model, cardinality $2^{\aleph_{0}}$ when using the
standard binary model, cardinality $\aleph_{0}!$ when using a factorial
numbering system and cardinality $\aleph_{0}^{\aleph_{0}}$ when using
the applicative numbering model.

With this new knowledge, the question now is how to find a proper
representation for the real numbers that is both valid and exhaustive.

\section{A proper model for the set of real numbers\label{sec:A-model-for-R}}

As mentioned previously, the standard positional fractional model
is not a proper model for $\mathbb{R}$, which means it can't be used
to well-order $\mathbb{R}$. However, it is still a valid model, so
it can be used to write down any individual real number we may need
to point to. Now, remember that the positional fractional notation
is just a model for the set of real numbers, not the set of real numbers
itself. In order to find a proper model for $\mathbb{R}$ we need
first to understand what $\mathbb{R}$ is.

The most compelling assumption is that $\mathbb{R}$ is the collection
of numbers (objects that satisfy the properties of the elements of
$\mathbb{R}$) produced by a set of mathematical algorithms. Huntington
in 1917 (\cite{Huntington2003}, p.16) writes ``it should be noticed
that what we are here required to grasp is not the infinite totality
of digits in the decimal fraction, but simply the rule by which those
digits are determined''--and those rules are the algorithms that
ultimately generate the numbers. By that time, however, there was
no general way to encode algorithms (programs), and Huntington contemporaries
tended to use ``...'' to indicate ``and so on'', assuming the
reader would be able to ``grasp'' the rule governing the omitted
digits of a number string. Yet without a very specific context the
use of the ``...'' notation is of course always ambiguous; but even
if it wasn't, strings of digits are not a good choice for encoding
algorithms.

The idea of numbers as the results of algorithms was further developed
by Turing in 1936, with the introduction of Turing machines as a general
way to encode the algorithms \cite{Turing1937}. A proper modern alternative
to model these ``rules'' that define the real numbers could be,
for instance, the set of strings that output a scalar real number
when parsed by Mathematica or Maple. This is of course a conjecture,
but note that it is not very different from the Church-Turing thesis,
which assumes without proof that any computable function is Turing-computable.
Under such an assumption we would be using a finite alphabet $\Sigma$
(consisting, for instance, of the ASCII character set) and restricting
ourselves to finite strings, since a program of infinite length never
yields a final return value. The cardinality of this algorithmic representation
$R_{alg}\left(\mathbb{R}\right)$ would therefore be $\aleph_{0}$.

Following this idea, notice that there is little difference between
the fractional notation string of a number \emph{x} and an imperative
program that implements its standard Val function. In fact, it can
be argued without much controversy that the string ``123.5'' is
but a shorthand for the algorithm ``\emph{x} = 1�100, \emph{x} =
\emph{x} + 2�10, \emph{x} = \emph{x} + 3�1, \emph{x} = \emph{x} +
5/10, return \emph{x}''.

Now, to see if the numbers-as-algorithms representation can constitute
a proper model $M_{alg}$ for $\mathbb{R}$, we need to check if it
is valid and exhaustive. Since the Val function for the standard positional
numbering model is but a simple algorithm, it is clear that the representation
of real numbers as the result of algorithms is valid. To see if this
representation is exhaustive we need to show that every real number
that may be the result of an infinite algorithm can be proved to be
equivalent to the result of another finite algorithm. This is easy
(at this point, after the big conjecture made earlier), because the
only way to determine if the result of a given infinite algorithm
is equivalent to a real number is because there is a finite process
we can use to verify the infinite algorithm converges to that particular
real number; that means there is an equivalent finite algorithm to
produce the number. Here we need to keep in mind that the infinite
algorithms we have to consider are not only infinite strings of decimal
digits but arbitrary infinite strings generated in the Kleene closure
of our alphabet $\Sigma$. This means that most of those infinite
strings will be malformed and will map to 0 or some predefined value;
others may diverge to infinity or include a division by zero or equivalently
invalid operations along their infinite list of instructions. The
important thing is that those infinite strings that actually represent
real numbers are those that can be reduced to equivalent finite strings
representing those same numbers; otherwise we would generate a contradiction:
there could be a specific infinite string that we know converges to
a real number but there is no finite algorithm (mathematical proof)
to demonstrate the equivalence.

\section{The Continuum Hypothesis}

The original form of the Continuum Hypothesis as posed by Cantor in
1878 was concerned about the possibility of other cardinals existing
between $\aleph_{0}$ and $\mathfrak{c}$. This is now known as the
Weak Continuum Hypothesis (WCH) \cite{Moore2011}. Since we have shown
that $\aleph_{0}=2^{\aleph_{0}}$ we can give a positive answer to
the WCH in the sense that there can't be any other cardinals in between,
because $\aleph_{0}$ and $\mathfrak{c}$ are equivalent.

As for what is now commonly called---after Hilbert---the Continuum
Hypothesis (CH), which questions if $\mathfrak{c}=\aleph_{1}$, what
we have seen so far reduces the CH to whether the statement $\aleph_{0}=\aleph_{1}$
is true or not. At this point we need to remember that transfinite
number theory was originally motivated by the apparent existence of
sets of higher cardinality than $\mathbb{N}$, but since we have shown
that the power-set operation for infinite sets does not actually produce
\textquotedblleft{}bigger\textquotedblright{} sets, there doesn\textquoteright{}t
seem to be any evidence or need for infinities \textquotedblleft{}bigger\textquotedblright{}
than $\aleph_{0}$. In any case, if we follow Cantor's definition
that $\aleph_{1}$ is the cardinality of the set of countable ordinals
$\omega_{1}$ (\cite{Cantor1915}, p.169;\cite{Huntington2003}, p.76)
then we just need to see if the cardinality of $\omega_{1}$ can be
proved equipollent to $\aleph_{0}$ or not.

\subsection{The cardinality of the set of ordinal numbers}

The ordinal numbers are an extension of the elements of $\mathbb{N}$
assuming the successor operator can continue after exhausting them.
For this reason, a new element $\omega$ is introduced. Combinations
involving natural numbers, $\omega$ and the operations +, $\times$
and \textasciicircum{} are formed until again we reach a new limiting
expression. At this point a new element $\varepsilon_{0}$ is introduced.
The same pattern is repeated indefinitely, adding new elements as
required. This necessarily results in a set of finite strings (formulas)
composed of symbols from an infinite alphabet $\Sigma$, which as
shown in section \ref{sec:Counter-Cantor-Theorem} yields a set of
cardinality $\aleph_{0}^{\aleph_{0}}$. A potential problem with this
approach could be the need to comprehend transfinitely many ordinal
numbers (see \cite{Conway1996}, p. 274), which would in principle
exhaust an alphabet of just $\aleph_{0}$ symbols; a straightforward
solution is to use a transfinite ordinal $\gamma$ to index a general
symbol $\varepsilon_{\gamma}$. A similar situation may happen with
the operators involved in the formulas, so we can collect them under
an indexed hyper-operator $\mathrm{H_{\gamma}}=\left\{ +,\,\times,\,\textrm{\textasciicircum},\,\uparrow\uparrow,\,\ldots\right\} $
(see \cite{Goodstein1947,Knuth1976}).

A model $M\left(\omega_{1}\right)$ for the set of ordinals $\omega_{1}$
thus requires an unending alphabet $\Sigma=\left\{ 1,\,\mathrm{H},\,(,\,),\,\omega,\,\varepsilon,\,\ldots\right\} $
with new symbols added as needed. Using unary notation (base \emph{b}=1)
any natural number can be represented as a concatenation of the symbol
1.

Note that here we are just concerned with whether any ordinal number
can be written in some way, independently of the underlying theory
\emph{T} being consistent or meaningful. With this in mind, we may
just assume that any new definition or constructive theorem added
to \emph{T} will result in the addition of a new symbol or group of
symbols to the alphabet $\Sigma$. As shown above, even the addition
of transfinitely many new objects can be done with just one new symbol
indexed by a transfinite ordinal. The assumption of $\Sigma$ having
$\aleph_{0}$ symbols simply means that we have at our disposal as
many symbols as we may need. In this sense, the ultimate question
regarding the cardinality of the set of ordinals is equivalent to
asking if we can always construct a finite string that unequivocally
refers to any specific ordinal knowing that we can use as many symbols
as we want.

Cantor established that any ordinal number can be represented using
a polynomial of ordinal numbers, known as the Cantor Normal Form (CNF)
of the number. The constructor function for the model thus simply
needs to write down the CNF of any given number using the symbols
in $\Sigma$; since there are many equivalent strings to any given
input we can select the one that is minimal by some predefined criterion
(e.g. coming first in the applicative ordering). Similarly the value
function for the model takes an input string of symbols from $\Sigma$
and outputs the equivalent CNF polynomial. In summary we have:

\begin{equation}
M\left(\omega_{1}\right)=\left\{ \begin{array}{l}
R\left(\omega_{1}\right)=\Sigma^{*};\thinspace\\
\mathrm{Rep}\left(x\in\omega_{1}\right)=\min\left\{ r\in R\,|\,\mathrm{Val}\left(r\right)=x\right\} ;\thinspace\\
\mathrm{Val}\left(r\in R\right)=\mathrm{CNF}\left(r\right);\thinspace\\
\mathrm{\Sigma=}\left\{ 1,\,\mathrm{H},\,(,\,),\,\omega,\,\varepsilon,\,\ldots\right\} 
\end{array}\right\} 
\end{equation}

Now we need to determine if such a model is valid and exhaustive.

One thing that is clear is that having as many symbols as we want
gives us a lot of freedom. It is difficult to conceive needing more
than an infinity of symbols, since that would mean the underlying
theory \emph{T} has more than infinitely many definitions or theorems,
something that is simply impossible to comprehend or conceive for
a human being.

The closure of the set of strings formed with the current set of symbols
necessarily yields a single maximal limit element (proof: the set
of ordinals is well ordered, so if there are two limiting expressions
one of them has to be maximal or otherwise both tend to the same limit).
A new single symbol is then introduced to represent that limit element
and the process continues indefinitely. This process of generating
the alphabet $\Sigma$, which adds a new element at a time, is analogous
to the generation of the natural numbers with the successor operation,
which necessarily yields a cardinality of $\aleph_{0}$ for $\Sigma$.

So, is this model valid? Well, for any input ordinal in CNF notation
there is a string in $M\left(\omega_{1}\right)$ with replicates the
CNF formula replacing its symbols with symbols from $\Sigma$. Also,
assuming there is an ordinal not covered by $M\left(\omega_{1}\right)$
we can just add another symbol to $\Sigma$; and note that even if
$\Sigma$ already has $\aleph_{0}$ elements we can still add another
denumerable infinity of symbols without problems by using the Hilbert's
Hotel strategy. The conclusion must therefore be that $M\left(\omega_{1}\right)$
is valid.

As for the model being exhaustive, in the boundary of the closure
of $R\left(\omega_{1}\right)$ there are limiting strings which are
either invalid (malformed strings, which are discarded) or valid,
in which case they must correspond to another simpler symbol in the
alphabet (because, as we just mentioned, that is specifically the
definition of that other symbol). We must again conclude that the
model is exhaustive.

Therefore, $M\left(\omega_{1}\right)$ is a proper model for $\omega_{1}$
and thus $R\left(\omega_{1}\right)$ and $\omega_{1}$ are equipollent,
which means they have the same cardinality. As a consequence $\aleph_{0}^{\aleph_{0}}=\aleph_{1}$
and, since we have proved in section \ref{sec:Counter-Cantor-Theorem}
that a set of cardinality $\aleph_{0}^{\aleph_{0}}$ can be put in
correspondence with the set of natural numbers, we must conclude $\aleph_{1}=\aleph_{0}$.

There is another way to see that the set of ordinals is of cardinality
$\aleph_{0}$ by showing that it is a denumerable collection of denumerable
sets: let $\mathcal{O}_{1}$ be the set of strings formed using only
the first symbol of $\Sigma$, $\mathcal{O}_{2}$ the set of strings
formed using the first two symbols, and so on; all the $\mathcal{O}_{i}$
sets are denumerable, and there are as many of them as symbols in
$\Sigma$ (which is denumerable), and since a union of denumerably
many denumerable sets is denumerable we have $\left|\omega_{1}\right|=\aleph_{0}$.
To make the $\mathcal{O}_{i}$ sets disjoint just consider the collection
$\mathcal{O}_{i}^{'}=\mathcal{O}_{i}\setminus\mathcal{O}_{i-1}$ instead,
which guarantees that $\mathcal{O}_{i}^{'}$ only has strings involving
the first \textit{i} symbols where the $i^{\mathrm{th}}$ symbol appears
at least once.

So, again, we see that there is no need or justification for introducing
a new cardinal for the ``size'' of the set of ordinals.

\bibliographystyle{ieeetr}
\bibliography{reference}

\end{document}